\documentclass[11pt]{article}
\usepackage{epic,latexsym,amssymb}
\usepackage{color}
\usepackage{tikz}
\usepackage{comment}
\usepackage{lineno}
\usepackage{pstricks,pst-node,pst-plot}
\usetikzlibrary{patterns}
\usepackage{geometry,graphicx}
\usepackage{todonotes}
\usepackage{subcaption}
\usepackage{hyperref}

\usetikzlibrary{calc}

\usepackage{amsmath,amsthm}

\newtheorem{theorem}{Theorem}
\newtheorem{lemma}[theorem]{Lemma}

\newtheorem{claim}[theorem]{Claim}

\newtheorem{observation}[theorem]{Observation}
\newtheorem{proposition}[theorem]{Proposition}

\newtheorem{rmk}[theorem]{Remark}
\newtheorem{definition}[theorem]{Definition}

\newcommand{\cT}{\mathcal{T}}
\newcommand{\bk}{{\bf{k}}}
\newcommand{\N}{\mathbb{N}}
\newcommand{\Nstar}{\mathbb{N}_\star^6}
\newcommand{\smallqed}{{\tiny ($\Box$)}}

\newcommand{\diam}{{\rm diam}}

\newcommand{\ioc}{\gamma^{{\rm {\small IOC}}}}

\newcommand{\proofclaim}{\noindent\textit{Proof of claim. }}
\newcommand{\2}{ \vspace{0.2cm} }
\newcommand{\1}{ \vspace{0.1cm} }

\let\oldenumerate\enumerate
\renewcommand{\enumerate}{
  \oldenumerate
  \setlength{\itemsep}{0pt}
  \setlength{\parskip}{0pt}
  \setlength{\parsep}{0pt}
}

\begin{document}

\title{Identifying open codes in trees and 4-cycle-free graphs\\of given maximum degree}

\author{Dipayan Chakraborty$^{a,b,c,d,e}$, Florent Foucaud$^{a}$, \, and \, Michael A. Henning$^{b}$\\ \\
$^a$ Universit\'{e} Clermont Auvergne, CNRS, Clermont Auvergne INP,\\
Mines Saint-\'Etienne, LIMOS, 63000 Clermont-Ferrand, France \\
\small \tt Email: dipayancha@gmail.com \\
\small \tt Email: florent.foucaud@uca.fr \\
\\
$^{b}$ Department of Mathematics and Applied Mathematics \\
University of Johannesburg, South Africa\\
\small \tt Email: mahenning@uj.ac.za \\
\\
$^{c}$ Department of Computer Science and Mathematics,\\
Lebanese American University, Beirut, Lebanon\\
\\
$^{d}$ Centrale M\'{e}diterran\'{e}e, Technopôle Château-Gombert,\\
13013 Marseille, France\\
\\
$^{e}$ Laboratoire d'Informatique et de Systèmes, UMR 7020, Marseille, France
}

\date{}
\maketitle

\begin{abstract}
An identifying open code of a graph $G$ is a set $S$ of vertices that is both a separating open code (that is, $N_G(u) \cap S \ne N_G(v) \cap S$ for all distinct vertices $u$ and $v$ in $G$) and a total dominating set (that is, $N(v) \cap S \ne \emptyset$ for all vertices~$v$ in $G$). Such a set exists if and only if the graph $G$ is open twin-free and isolate-free; and the minimum cardinality of an identifying open code in an open twin-free and isolate-free graph $G$, its identifying open code number, is denoted by $\ioc(G)$.

We study the identifying open code number of a graph, in relation with its order and its maximum degree. For $\Delta$ a fixed integer at least~$3$, we show that if $G$ is a connected graph of order~$n \ge 5$ that contains no $4$-cycle and is open twin-free with maximum degree bounded above by~$\Delta$, then  $\ioc(G) \le \left( \frac{2\Delta - 1}{2\Delta} \right) n$, unless $G$ is obtained from a star $K_{1,\Delta}$ by subdividing every edge exactly once. Moreover, we show that the bound is best possible by constructing graphs that reach the bound when $\Delta=3$, and nearly best possible by another construction when $\Delta\geq 4$, with identifying open code numbers $\left(\frac{2\Delta - 4}{2\Delta-3} \right) n$.
\end{abstract}


\section{Introduction}
\label{S:into}

The goal of this paper is to investigate the possible size of a smallest identifying open code in a graph of given maximum degree. Identifying open codes in graphs define a concept that is simultaneously part of two rich areas in discrete mathematics: domination problems and identification problems.

\subsection{Definitions}

A set $S$ of vertices in a graph $G$ is a \emph{dominating set} if every vertex not in $S$ is adjacent to a vertex in~$S$. Further if every vertex in $G$ is adjacent to some other vertex in $S$, then $S$ is a \emph{total dominating set}, abbreviated TD-set of $G$. The \emph{total domination number} $\gamma_t(G)$ of $G$ is the minimum cardinality of a TD-set of $G$. A TD-set of cardinality $\gamma_t(G)$ is called a $\gamma_t$-\emph{set of $G$}. A vertex $v$ \emph{totally dominates} another vertex $u$ if they are adjacent vertices. More generally, if $X$ and $Y$ are subsets of vertices in $G$, then the set $X$ \emph{totally dominates} the set $Y$ in $G$ if every vertex in $Y$ is adjacent to at least one vertex in $X$. In particular, if $X$ totally dominates $V(G)$, then $X$ is a TD-set of $G$. For recent books on domination in graphs, we refer the reader to~\cite{HaHeHe-20,HaHeHe-21,HaHeHe-23,HeYe-book}.

For graph theory notation and terminology, we generally follow~\cite{HeYe-book}.  Specifically, let $G$ be a graph with vertex set $V(G)$ and edge set $E(G)$. The order and size of a graph $G$, denoted $n(G)$ and $m(G)$, is $|V(G)|$ and  $|E(G)|$, respectively. Two vertices $u$ and $v$ of $G$ are \emph{adjacent} if $uv \in E(G)$, and are called \emph{neighbors}. The \emph{open neighborhood} $N_G(v)$ of a vertex $v$ in $G$ is the set of neighbors of $v$, while the \emph{closed neighborhood} of $v$ is the set $N_G[v] = \{v\} \cup N_G(v)$. The \emph{degree} of a vertex $v$ in $G$ is the number of neighbors $v$ in $G$, and is denoted by $\deg_G(v)$, and so $\deg_G(v) = |N_G(v)|$. An \emph{isolated vertex} in $G$ is a vertex of degree~$0$. A graph is \emph{isolate}-\emph{free} if it contains no isolated vertices. A \emph{leaf} in $G$ is a vertex of degree~$1$ in $G$, and a \emph{support vertex} in $G$ is a vertex with at least one leaf neighbor. A \emph{strong support vertex} in $G$ is a vertex with at least two leaf neighbors. The minimum and maximum degrees in $G$ are denoted by $\delta(G)$ and $\Delta(G)$, respectively.

A \emph{cycle} on $n$ vertices is denoted by $C_n$ and a \emph{path} on $n$ vertices by $P_n$. The  \emph{complete graph} on $n$ vertices is denoted by $K_n$, while the \emph{complete bipartite graph} with one partite set of size~$n$ and the other of size~$m$ is denoted by $K_{n,m}$. A \emph{star} is a complete bipartite graph $K_{1,n}$. A \emph{cubic graph} is graph in which every vertex has degree~$3$, while a \emph{subcubic graph} is a graph with maximum degree at most~$3$. A \emph{subcubic tree} is a tree with maximum degree at most~$3$. The \emph{diameter} of a graph $G$, denoted $\diam(G)$,  is the maximum distance among all pairs of vertices of $G$. In particular, the diameter of a tree $T$ is the length of a longest path in $T$.

A \emph{rooted tree} $T$ distinguishes one vertex $r$ called the \emph{root}. Let $T$ be a tree rooted at vertex~$r$.  For each vertex  $v \ne r$ of $T$, the \emph{parent} of $v$ is the neighbor of $v$ on the unique $(r,v)$-path, while a \emph{child} of $v$ is any other neighbor of $v$. The root $r$ does not have a parent in $T$ and all its neighbors are its children.  A \emph{descendant} of $v$ is a vertex $x$ such that the unique $(r,x)$-path contains $v$. Thus, every child of $v$ is a descendant of $v$. A \emph{grandchild} of $v$ is a descendant of $v$ at distance~$2$ from~$v$. Let $C(v)$ and $D(v)$ denote the set of children and descendants, respectively, of $v$, and we define $D[v] = D(v) \cup \{v\}$. The \emph{maximal subtree} at $v$, denoted $T_v$, is the subtree of $T$ induced by the set $D[v]$. For $k \ge 1$ an integer, we let $[k]$ denote the set $\{1,\ldots,k\}$ and we let $[k]_0 = [k] \cup \{0\} = \{0,1,\ldots,k\}$.

A graph $G$ has a TD-set if and only if it is isolate-free. Hence throughout this paper all graphs are isolate-free, unless otherwise stated. A \emph{separating open code} in $G$ is a set $S$ of vertices such that $N_G(u) \cap S \ne N_G(v) \cap S$ for all distinct vertices $u$ and $v$ in $G$. An \emph{identifying open code} of $G$, abbreviated IO-code, is a set $S$ that is both a separating open code and a TD-set. Two vertices in $G$ are \emph{open twins} if they have the same open neighborhood. A graph is \emph{open twin-free} if it has no open twins.  An isolate-free graph has an IO-code if and only if it is open twin-free. The minimum cardinality of an IO-code in an open twin-free graph $G$ is denoted by $\ioc(G)$ and called the \emph{identifying open code number} of $G$ (IO-code number for short).

\subsection{Previous work}

An IO-code has also been called an \emph{open (neighborhood) locating-dominating set}~\cite{SeSl10,SeSl11} or an \emph{identifying code with nontransmitting faulty vertices}~\cite{HoLaRa02} in the literature. We follow the terminology from~\cite{HeYe-14}. The problem of ``identifying open codes'' was introduced by Honkala, Laihonen and Ranto~\cite{HoLaRa02} in the context of coding theory for binary hypercubes and further studied, for example, in~\cite{ArBiLuWa-20,ArBiLuWa-22,CaCoEFo-22,ChFoPaWa-24,Che-14,FoGhSh-21,FoGhSh-24,FoMeNaPaVa17a,FoMeNaPaVa17b,Giv-22,HeYe-14,JeSe-23,Kin-15,Pan-17,SeSl10,SeSl11}. This concept is closely related to the extensively studied concepts of identifying codes and locating-dominating sets in graphs, and is part of a vast literature on general identification problems in discrete structures, with many theoretical and practical applications. For a survey on (open) identifying codes and locating (total) domination in graphs we refer the reader to the book chapter by Lobstein, Hudry, and Charon~\cite{LoHuCh-20}, and for an online bibliography maintained by Jean and Lobstein on these topics, see~\cite{biblio}.

Regarding upper bounds on the IO-code number, the graphs $G$ of order $n$ with $\ioc(G)=n$ have been characterized by Foucaud, Ghareghani, Roshany-Tabrizi and Sharifani in~\cite{FoGhSh-21} as the family of \emph{half-graphs} defined in~\cite{EH84} (the half-graph $H_k$ is the bipartite graph with both parts of size $k$, where each part can be ordered so that the open neighborhoods of consecutive vertices differ by exactly one vertex). In~\cite{SeSl11}, Seo and Slater characterized the trees $T$ of order $n$ with $\ioc(T)=n-1$.

Such upper bounds have also been studied for classic identifying codes, where the open neighborhood is replaced by the closed neighourhood in their definition. In particular, some works have explored the best possible upper bound for the optimal size of an identifying code, depending on the maximum degree and the order of the graph~\cite{ChFoLa23,ChFoHeLa-24a,ChFoHeLa-24b,FoKla-12,FoPe-12}. Those works have inspired the current article.

Further motivation comes from the work by Henning and Yeo~\cite{HeYe-14}, who initiated the study of upper bounds on the IO-code number of graphs of given maximum degree, focusing on regular graphs. They proved, using an interplay with transversal of hypegraphs, that $\ioc(G)\le \frac{3}{4}n$ for every open twin-free connected cubic graph $G$ and this bound is tight for the complete graph $K_4$ and the hypercube of dimension~$3$. They also suggested that perhaps, more generally,  $\ioc(G)\le \left(\frac{\Delta}{\Delta+1}\right)n$ holds for every connected open twin-free $\Delta$-regular graph $G$ of order $n$. Our goal is to investigate what happens when the graph is non-regular. We will see that, interestingly, the above conjectured bound for the IO-code number of regular graphs must be modified for non-regular graphs.

\subsection{Main result}

Let $G$ be a connected graph of order~$n \ge 2$ that is open twin-free. If $n = 2$, then $G \cong P_2$ and $\ioc(G) = 2 = n$. If $n = 3$, then $G \cong K_3$ and $\ioc(G) = 2 = \frac{2}{3}n$. If $n = 4$, then $G \cong P_4$ and $\ioc(G) = 4 = n$. Our aim is to obtain a best possible upper bound on IO-codes in open twin-free connected graphs $G$ of order~$n \ge 5$ that contains no $4$-cycles and have bounded maximum degree.

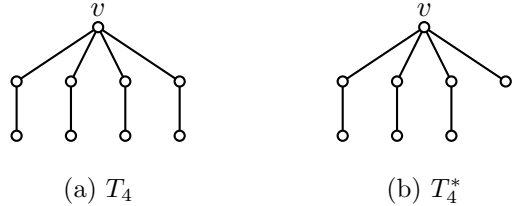
\begin{figure}[htb]
\begin{center}
\begin{tikzpicture}[scale=.85,style=thick,x=0.85cm,y=0.85cm]
\def\vr{2.25pt}
\path (0,0) coordinate (u1);
\path (0,1) coordinate (v1);
\path (1,0) coordinate (u2);
\path (1,1) coordinate (v2);
\path (2,0) coordinate (u3);
\path (2,1) coordinate (v3);
\path (3,0) coordinate (u4);
\path (3,1) coordinate (v4);
\path (1.5,2) coordinate (v);
\draw (v)--(v1)--(u1);
\draw (v)--(v2)--(u2);
\draw (v)--(v3)--(u3);
\draw (v)--(v4)--(u4);
%
\draw (v1) [fill=white] circle (\vr);
\draw (v2) [fill=white] circle (\vr);
\draw (v3) [fill=white] circle (\vr);
\draw (v4) [fill=white] circle (\vr);
\draw (u1) [fill=white] circle (\vr);
\draw (u2) [fill=white] circle (\vr);
\draw (u3) [fill=white] circle (\vr);
\draw (u4) [fill=white] circle (\vr);
\draw (v) [fill=white] circle (\vr);
\draw[anchor = south] (v) node {$v$};
\draw (1.5,-1) node {{\small (a) $T_4$}};
\path (6,0) coordinate (u1);
\path (6,1) coordinate (v1);
\path (7,0) coordinate (u2);
\path (7,1) coordinate (v2);
\path (8,0) coordinate (u3);
\path (8,1) coordinate (v3);
\path (9,1) coordinate (v4);
\path (7.5,2) coordinate (v);
\draw (v)--(v1)--(u1);
\draw (v)--(v2)--(u2);
\draw (v)--(v3)--(u3);
\draw (v)--(v4);
%
\draw (v1) [fill=white] circle (\vr);
\draw (v2) [fill=white] circle (\vr);
\draw (v3) [fill=white] circle (\vr);
\draw (v4) [fill=white] circle (\vr);
\draw (u1) [fill=white] circle (\vr);
\draw (u2) [fill=white] circle (\vr);
\draw (u3) [fill=white] circle (\vr);
\draw (v) [fill=white] circle (\vr);
\draw[anchor = south] (v) node {$v$};
\draw (7.5,-1) node {{\small (b) $T_4^*$}};
\end{tikzpicture}
\caption{The subdivided star $T_4 = S(K_{1,4})$.}
\label{fig:T4}
\end{center}
\end{figure}

For $\Delta \ge 3$, a \emph{subdivided star} $T_{\Delta} = S(K_{1,\Delta})$ is the tree obtained from a star $K_{1,\Delta}$ with universal vertex $v$, say, by subdividing every edge exactly once. The resulting tree $T_{\Delta}$ has order~$n = 2\Delta +1$ and maximum degree $\Delta(T_{\Delta}) = \Delta$. For example, when $\Delta = 4$ the subdivided star $T_4 = S(K_{1,4})$ is illustrated in Figure~\ref{fig:T4}(a).
For $\Delta \ge 3$, a \emph{reduced subdivided star} $T_{\Delta}^*$ is the tree obtained from a subdivided star $K_{1,\Delta}$ by removing exactly one leaf. The resulting tree $T_{\Delta}^*$ has order~$n = 2\Delta$ and maximum degree $\Delta(T_{\Delta}) = \Delta$. For example, when $\Delta = 4$ the reduced subdivided star $T_4^*$ is illustrated in Figure~\ref{fig:T4}(b). In both the subdivided and the reduced subdivided tree, the vertex $v$ is called the \emph{central vertex} of the tree.

Our main result is to prove that, not only for trees but also for graphs with no $4$-cycles, among all graphs of given maximum degree, the subdivided stars described above require the largest fraction of their vertex set in any optimal IO-code. More precisely, we prove the following statement.

\begin{theorem}
\label{thm:main1}
For $\Delta \ge 3$ a fixed integer, if $G \not \cong T_\Delta$ is an open twin-free connected graph of order~$n \ge 5$ that contains no $4$-cycles and satisfies $\Delta(G) \le \Delta$, then
\[
\ioc(G) \le \left( \frac{2\Delta - 1}{2\Delta} \right) n.
\]
On the other hand, when $G \cong  T_{\Delta}$, then we have $\ioc(G) = \left( \frac{2\Delta}{2\Delta + 1} \right)n$.
\end{theorem}

Moreover, we show that the above bound is nearly tight for every value of $\Delta\ge 3$ by providing a construction of graphs with IO-code numbers equal to $\left(\frac{2\Delta-4}{2\Delta-3}\right)n$. Furthermore, when $\Delta=3$, we give another construction that provides infinitely many connected graphs for which the bound is attained.

\subsection{Organization of the paper}

We first prove Theorem~\ref{thm:main1} for trees in Section~\ref{S:tree-proof}: this is in fact the essential part of the proof. We then prove Theorem~\ref{thm:main1} in full generality in Section~\ref{S:proof}. We provide constructions to show the tightness of the bound in Section~\ref{S:construct}, and we conclude in Section~\ref{S:conclude} with directions for future research. 

\section{Trees}
\label{S:tree-proof}

In this section, we prove Theorem~\ref{thm:main1} for trees. Towards this goal, we first define some special families of trees that will be essential in the proof.

\subsection{A special family $\cT$ of trees}

In this section, we define a special family of trees that we will often need when proving our main result. In order to define our special family of trees, we introduce some additional notations. Let $r$ be a specified vertex in a tree $T$. We define next several types of attachments at the vertex~$r$ that we use to build larger trees. In all cases, we call the vertex of the attachment that is joined to $r$ the \emph{link vertex} of the attachment. \\[-20pt]
\begin{enumerate}
\item[$\bullet$] For $i \in [4]$, an \emph{attachment of Type-$i$ at $r$} is an operation that adds a path $P_{i}$ to $T$ and joins one of its leaves to $r$. \1
\item[$\bullet$] An \emph{attachment of Type-$5$ at $r$} is an operation that adds a path $P_4$ to $T$ and joins one of its support vertices to $r$. \1
\item[$\bullet$] An \emph{attachment of Type-$6$ at $r$} is an operation that adds a star $K_{1,3}$ with one edge subdivided and joins one of its leaves that is adjacent to the vertex of degree~$3$ to $r$.
\end{enumerate}
Let $\N$ be the set of all non-negative integers and let
\begin{align} \label{eq_Nstar}
\begin{split}
\Nstar = \{ &\bk = (k_1,k_2,k_3,k_4,k_5,k_6) \in \N^6 : k_1 \in \{0,1\} \text{ and }\\ &\bk \ne  (1,0,0,0,0,0), (0,1,0,0,0,0), (0,0,1,0,0,0), (1,1,0,0,0,0) \}.
\end{split}
\end{align}

From now on, let $\bk$ denote a vector in $\Nstar$, let $k_i$ denote the $i$th coordinate of $\bk$ and let $k = \sum_{i=1}^6 k_i$. We define $T(r;\bk)$ to be the tree obtained from a trivial tree $K_1$ whose vertex is named~$r$ by applying $k_i$ attachments of Type-$i$ at vertex $r$ for each $i \in [6]$. Notice that $\deg_T (r) = k$. We now define the set $\cT = \{T(r;\bk) : \bk \in \Nstar\}$ to be our special family of trees. We note that the forbidden configurations of the vector $\bk$ in the definition of $\Nstar$ exclude the trees $P_2$, $P_3$ and $P_4$ from the family $\cT$, as in our main result, we only consider graphs of order at least~$5$.

\begin{figure}[t]
\centering
\begin{subfigure}[t]{1\textwidth}
\centering
\begin{tikzpicture}[scale=.85,style=thick,x=0.85cm,y=0.85cm]
\def\vr{2.25pt}
\path (-0.5,3) coordinate (u0);
\path (1,3) coordinate (u1);
\path (1,2) coordinate (v1);
\path (2,3) coordinate (u2);
\path (2,2) coordinate (v2);
\path (3,3) coordinate (u3);
\path (3,2) coordinate (v3);
\path (4,3) coordinate (u4);
\path (4,2) coordinate (v4);
\path (4,1) coordinate (w4);
\path (5,3) coordinate (u5);
\path (5,2) coordinate (v5);
\path (5,1) coordinate (w5);
\path (6,3) coordinate (u8);
\path (6,2) coordinate (v8);
\path (6,1) coordinate (w8);
\path (6,0) coordinate (x8);
\path (7,3) coordinate (u9);
\path (7,2) coordinate (v9);
\path (7,1) coordinate (w9);
\path (7,0) coordinate (x9);
\path (8,3) coordinate (u10);
\path (8,2) coordinate (v10);
\path (8,1) coordinate (w10);
\path (8,0) coordinate (x10);
\path (9,3) coordinate (u6);
\path (9,2) coordinate (v6);
\path (9.5,2) coordinate (v6p);
\path (9,1) coordinate (w6);
\path (10.5,3) coordinate (u7);
\path (10.5,2) coordinate (v7);
\path (11,2) coordinate (v7p);
\path (10.5,1) coordinate (w7);
\path (12.5,3) coordinate (u11);
\path (12.5,2) coordinate (v11);
\path (12.5,1) coordinate (w11);
\path (13,1) coordinate (w11p);
\path (12.5,0) coordinate (x11);
\path (14,3) coordinate (u12);
\path (14,2) coordinate (v12);
\path (14,1) coordinate (w12);
\path (14.5,1) coordinate (w12p);
\path (14,0) coordinate (x12);
\path (6.5,5.5) coordinate (v);
\draw (v)--(u0);
\draw (v)--(u1)--(v1);
\draw (v)--(u2)--(v2);
\draw (v)--(u3)--(v3);
\draw (v)--(u4)--(v4)--(w4);
\draw (v)--(u5)--(v5)--(w5);
\draw (v)--(u6)--(v6)--(w6);
\draw (u6)--(v6p);
\draw (v)--(u7)--(v7)--(w7);
\draw (u7)--(v7p);
\draw (v)--(u8)--(v8)--(w8)--(x8);
\draw (v)--(u9)--(v9)--(w9)--(x9);
\draw (v)--(u10)--(v10)--(w10)--(x10);
\draw (v)--(u11)--(v11)--(w11)--(x11);
\draw (v11)--(w11p);
\draw (v)--(u12)--(v12)--(w12)--(x12);
\draw (v12)--(w12p);
\draw (u0) [fill=white] circle (\vr);
\draw (u1) [fill=black] circle (\vr);
\draw (u2) [fill=black] circle (\vr);
\draw (u3) [fill=black] circle (\vr);
\draw (u4) [fill=black] circle (\vr);
\draw (u5) [fill=black] circle (\vr);
\draw (u6) [fill=black] circle (\vr);
\draw (u7) [fill=black] circle (\vr);
\draw (u8) [fill=black] circle (\vr);
\draw (u9) [fill=black] circle (\vr);
\draw (u10) [fill=black] circle (\vr);
\draw (u11) [fill=black] circle (\vr);
\draw (u12) [fill=black] circle (\vr);
\draw (v) [fill=black] circle (\vr);
\draw (v1) [fill=black] circle (\vr);
\draw (v2) [fill=black] circle (\vr);
\draw (v3) [fill=black] circle (\vr);
\draw (v4) [fill=black] circle (\vr);
\draw (v5) [fill=black] circle (\vr);
\draw (v6) [fill=black] circle (\vr);
\draw (v6p) [fill=white] circle (\vr);
\draw (v7) [fill=black] circle (\vr);
\draw (v7p) [fill=white] circle (\vr);
\draw (v8) [fill=black] circle (\vr);
\draw (v9) [fill=black] circle (\vr);
\draw (v10) [fill=black] circle (\vr);
\draw (v11) [fill=black] circle (\vr);
\draw (v12) [fill=black] circle (\vr);
\draw (w4) [fill=white] circle (\vr);
\draw (w5) [fill=white] circle (\vr);
\draw (w6) [fill=black] circle (\vr);
\draw (w7) [fill=black] circle (\vr);
\draw (w8) [fill=black] circle (\vr);
\draw (w9) [fill=black] circle (\vr);
\draw (w10) [fill=black] circle (\vr);
\draw (w11) [fill=black] circle (\vr);
\draw (w11p) [fill=white] circle (\vr);
\draw (w12) [fill=black] circle (\vr);
\draw (w12p) [fill=white] circle (\vr);
\draw (x8) [fill=white] circle (\vr);
\draw (x9) [fill=white] circle (\vr);
\draw (x10) [fill=white] circle (\vr);
\draw (x11) [fill=black] circle (\vr);
\draw (x12) [fill=black] circle (\vr);
\draw[anchor = south] (v) node {$r$};
\draw (-0.5,2.75) node {$\underbrace{\phantom{1}}$};
\draw (-0.5,2) node {$k_1 = 1$};
\draw (2,1.75) node {$\underbrace{\phantom{11111111}}$};
\draw (2,1) node {$k_2 = 3$};
\draw (4.5,0.75) node {$\underbrace{\phantom{11111}}$};
\draw (4.5,0) node {$k_3 = 2$};
\draw (7,-0.25) node {$\underbrace{\phantom{11111111}}$};
\draw (7,-1) node {$k_4 = 3$};
\draw (10,0.75) node {$\underbrace{\phantom{11111111}}$};
\draw (10,0) node {$k_5 = 2$};
\draw (13.25,-0.25) node {$\underbrace{\phantom{11111111}}$};
\draw (13.25,-1) node {$k_6 = 2$};
\end{tikzpicture}
\caption{The tree $T(r;1,3,2,3,2,2)$.}
\end{subfigure}
\vskip 0.25 cm 
\begin{subfigure}[t]{1\textwidth}
\centering
\begin{tikzpicture}[scale=.85,style=thick,x=0.85cm,y=0.85cm]
\def\vr{2.25pt}
\path (0,2) coordinate (u);
\path (1,3) coordinate (v);
\path (2,0) coordinate (v1);
\path (2,1) coordinate (v2);
\path (2,2) coordinate (v3);
\draw (v)--(v3)--(v2)--(v1);
\draw (v)--(u);
%
\draw (v1) [fill=white] circle (\vr);
\draw (v2) [fill=black] circle (\vr);
\draw (v3) [fill=black] circle (\vr);
\draw (v) [fill=black] circle (\vr);
\draw (u) [fill=black] circle (\vr);
\draw[anchor = south] (v) node {$r$};
\draw[anchor = north] (u) node {$u_1$};
\draw (0.8,-1) node {{\small (b) $T(r;1,0,1,0,0,0) \cong T_2$.}};
\path (6,2) coordinate (u);
\path (7,3) coordinate (v);
\path (8,0) coordinate (v1);
\path (8,1) coordinate (v2);
\path (9,1) coordinate (u2);
\path (8,2) coordinate (v3);
\draw (v)--(v3)--(v2)--(v1);
\draw (v)--(u);
\draw (v3)--(u2);
%
\draw (v1) [fill=black] circle (\vr);
\draw (v2) [fill=black] circle (\vr);
\draw (v3) [fill=black] circle (\vr);
\draw (v) [fill=black] circle (\vr);
\draw (u) [fill=black] circle (\vr);
\draw (u2) [fill=white] circle (\vr);
\draw[anchor = south] (v) node {$r$};
\draw[anchor = north] (u) node {$u_1$};
\draw (7.5,-1) node {{\small (c) $T(r;1,0,0,0,1,0) \cong T^*_3$.}};
\end{tikzpicture}
\caption*{}
\end{subfigure}
\caption{The canonical sets of a general tree and the trees $T_2$ and $T^*_3$ in $\cT$.}
\label{fig:special_trees}
\end{figure}
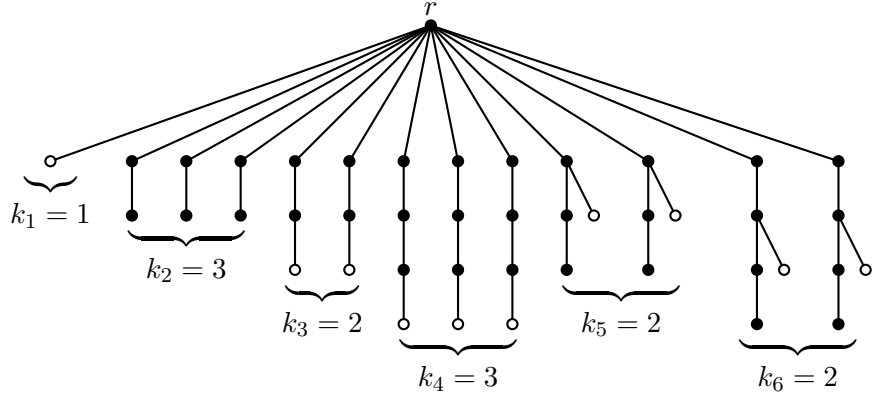
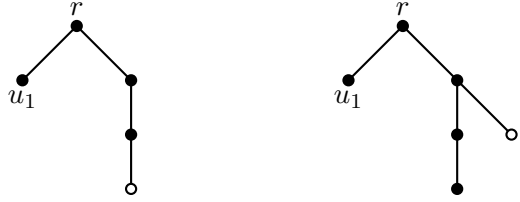

The tree $T(r;1,3,2,3,2,2)$, for example, is illustrated in Figure~\ref{fig:special_trees}(a). For $k \ge 2$, the subdivided star $T_k$ is isomorphic to the tree $T(r;0,k,0,0,0,0) \in \cT$ and also to the tree $T(r;1,0,1,0,0,0) \in \cT$ for $k = 2$ (see Figure~\ref{fig:special_trees}(c)). For $k \ge 3$, the reduced subdivided star $T^*_k$ is isomorphic to the tree $T(r;1,k-1,0,0,0,0) \in \cT$ and also to the tree $T(r;1,0,0,0,1,0)$ for $k=2$ (see Figure~\ref{fig:special_trees}(c)).

\begin{definition}
\label{defn1}
{\rm
We define the \emph{canonical set} $C$ of a tree $T = T(r;\bk) \in \cT$ as follows. For $\bk \notin \{(1,0,1,0,0,0,), (1,0,0,0,1,0)\}$, that is, when $T$ is neither $T_2$ nor $T^*_3$ with $r$ being a degree~$2$ support vertex of $T$ (as illustrated in Figures~\ref{fig:special_trees}(b) and (c), respectively), we define the canonical set $C$ of $T$ as follows. \\[-20pt]
\begin{enumerate}
\item[$\bullet$] Initially, we add the vertex~$r$ to the set $C$. \1
\item[$\bullet$] If $k_1 = 1$ and $k_2 \ge 1$, then we add to $C$ all vertices that belong to attachments of Type-$2$ at $r$, but we do not add to $C$ the vertex that belongs to the attachment of Type-$1$ at $r$.
\item[$\bullet$] If $k_1 = 0$ and $k_2 \ge 1$, then we add to $C$ all vertices that belong to attachments of Type-$2$ at $r$, except for one leaf at distance~$2$ from~$r$ in $T$ in exactly one attachment of Type-$2$.
\item[$\bullet$] We add to $C$ all vertices that belong to attachments of Type-$3$ at $r$ that are not leaves at distance~$3$ from~$r$ in $T$.
\item[$\bullet$] We add to $C$ all vertices that belong to attachments of Type-$4$ at $r$ that are not leaves at distance~$4$ from~$r$ in $T$.
\item[$\bullet$] We add to $C$ all vertices of attachments of Type-$5$ at $r$ except for the leaves at distance~$2$ from~$r$ in $T$.
\item[$\bullet$] We add to $C$ all vertices of attachments of Type-$6$ at $r$ except for the leaves at distance~$3$ from~$r$ in $T$.
\end{enumerate}

For $\bk \in \{(1,0,1,0,0,0), (1,0,0,0,1,0)\}$, that is, when $T$ is either $T_2$ or $T^*_3$ with $r$ being a degree~$2$ support vertex of $T$ (as in Figures~\ref{fig:special_trees}(b) and (c), respectively), we define the canonical set $C$ of $T$ as follows. \\[-20pt]
\begin{enumerate}
\item[$\bullet$] For $\bk = (1,0,1,0,0,0)$, add to $C$ all vertices of $T$ except for the leaf of $T$ at distance~$3$ of~$r$.
\item[$\bullet$] For $\bk = (1,0,1,0,0,0)$, add to $C$ all vertices of $T$ except for the leaf of $T$ at distance~$2$ of~$r$.
\end{enumerate}

For notational convenience, we shall abbreviate from now on the notations for the trees of the form $T(r;(1,0,1,0,0,0))$ as $T(r;1,0,1,0,0,0)$. For the trees $T(r;1,3,2,3,2,2)$, $T(r;1,0,1,0,0,0)$ and $T(r;1,0,0,0,1,0)$, for example, the shaded vertices in Figures~\ref{fig:special_trees}(a), \ref{fig:special_trees}(b) and \ref{fig:special_trees}(c), respectively, indicate the canonical set $C$ of the tree. Notice that the two trees $T(r;\bk)$ for $\bk \in \{(0,2,0,0,0,0), (1,0,1,0,0,0)\}$ (both isomorphic to $T_2$) have the same canonical sets as prescribed in Definition~\ref{defn1}. This property also holds for the trees $T(r;\bk)$ for $\bk \in \{(1,2,0,0,0,0), (1,0,0,0,1,0)\}$ (both isomorphic to $T^*_3$).
}
\end{definition}

\begin{observation}
\label{obser1}
The canonical set of any tree $T = T(r;\bk) \in \cT$ is also an IO-code of $T$.
\end{observation}

We next determine the IO-code numbers of (reduced) subdivided stars $T_{k}$ and~$T_{k}^*$.

\begin{proposition}
\label{prop:prop1}
For $k \ge 2$ a fixed integer and a tree $T$ of order~$n$, the following properties hold. \\ [-22pt]
\begin{enumerate}
\item[{\rm (a)}] If $T \cong T_{k}$, then $\ioc(T) = \left( \frac{2k}{2k + 1} \right)n$. \1
\item[{\rm (b)}] If $k \ge 3$ and $T \cong T_{k}^*$, then $\ioc(T) = \left( \frac{2k-1}{2k} \right)n$.
\end{enumerate}
Moreover, in both cases, the canonical set of $T$ is an optimal IO-code of $T$.
\end{proposition}
\begin{proof} (a) Let $T \cong T_{k}$ have order~$n$ where $k \ge 2$, and let $v$ be the central vertex of $T$ (of degree~$k$). Thus, $n = 2k + 1$. Let $S$ be an IO-code of $T$. Since $S$ is a TD-set of $T$, in order to totally dominate all the leaves in $T$ the set $S$ contains all support vertices of $T$. Since $S$ is a separating open code, in order to identify the support vertices of $T$ the set $S$ contains all, except possibly one leaf. If $S$ contains all leaves of $T$, then $|S| \ge 2k$. Suppose that $S$ does not contain all leaves of $T$. By our earlier observations, in this case there is exactly one leaf, $u$ say, that does not belong to $S$. Let $w$ be the support vertex adjacent to~$u$. Thus, $N_T(w) = \{u,v\}$. In order to totally dominate the vertex~$w$, we infer that $v \in S$, implying once again that $|S| \ge 2k$. Since $S$ is an arbitrary IO-code of $T$, this implies that $\ioc(T) \ge 2k$. Now, assuming $T \cong T(r;0,k,0,0,0,0) \in \cT$, let $C$ be the canonical set of the tree $T$ as in Definition~\ref{defn1}, that is, $C$ is the set of all vertices of $T$ except for one leaf of $T$. Then, by Observation~\ref{obser1}, the set $C$ is an IO-code and it satisfies $\ioc(T) \le |C| = 2k$. Consequently, $\ioc(T) = 2k = \left( \frac{2k}{2k + 1} \right)n$.

(b) Let $T \cong T_{k}^*$ have order~$n$ where $k \ge 3$, and let $v$ be the central vertex of $T$ (of degree~$k$) and let $x$ denote the leaf neighbor of $v$. Thus, $n = 2k$. Let $S$ be an IO-code of $T$. Since $S$ is a TD-set of $T$, the set $S$ contains all support vertices of $T$. Since $S$ is an IO-code of $T$, in order to identify the leaf $x$ from all other neighbors of $v$, the set $S$ contains all leaves of $T$ at distance~$2$ from~$v$. Thus, $S$ contains all vertices of $T$, except possibly for the leaf~$x$ of $T$. Thus, $|S| \ge 2k - 1$. Since $S$ is an arbitrary IO-code of $T$, this implies that $\ioc(T) \ge 2k - 1$. On the other hand, assuming $T \cong T(r;1,k-1,0,0,0,0)$, let $C$ be the canonical set of $T$ as in Definition~\ref{defn1}, that is, $C$ consists of all vertices of $T$ except for the leaf neighbor $x$ of $v$. Again, by Observation~\ref{obser1}, the set $C$ is an IO-code of $T$. Moreover, we have $\ioc(T) \le |C| = 2k - 1$. Consequently, $\ioc(T) = 2k - 1 = \left( \frac{2k - 1}{2k} \right)n$.
\end{proof}

\medskip 
We next prove that the bound of Theorem~\ref{thm:main1} holds for the trees in $\mathcal T$.

\begin{proposition}
\label{prop2_new}
Let $\Delta \ge 3$ be a fixed integer and let $T \cong T(r; \bf k) \in \cT$ be a tree of order $n \ge 5$, maximum degree $\Delta(T) \le \Delta$ such that $\deg_T(r) = k \ge 2$. If $T$ is not isomorphic to the subdivided star $T_\Delta$, then 
\[
\ioc(T) \le \left( \frac{2\Delta - 1}{2\Delta} \right) n.
\]
\end{proposition}
\begin{proof}
By supposition, $\Delta \ge \Delta(T) \ge \deg_T(r) = k$. Therefore, if $T \cong T_k$ for some $k \le \Delta - 1$, then the result follows from Proposition~\ref{prop:prop1}. Notice that $T$ cannot be isomorphic to the reduced subdivided star $T^*_2 \cong P_4$, since $T \in \cT$ and the latter does not contain the path $P_4$. On the other hand, if $T$ is isomorphic to the reduced subdivided star $T_k^*$ for $k \ge 3$, then again the desired upper bound follows from Proposition~\ref{prop:prop1}. Hence, we may assume that the tree $T$ is neither isomorphic to $T_k$ nor to $T^*_k$ for all $2 \le k \le \Delta$. Therefore, in particular, $T$ is neither the tree $T(r;1,0,1,0,0,0)$ as in Figure~\ref{fig:special_trees}(b) nor the tree $T(r;1,0,0,0,1,0)$ as in Figure~\ref{fig:special_trees}(c).

Let $C$ be the canonical set of $T \cong T(r;\bk)$ as in Definition~\ref{defn1}. Let $A$ be the set of vertices of $T$ consisting of the vertex~$r$ and all vertices that belong to attachments of Type-$1$ or of Type-$2$ at $r$, and let $B = V(T) \setminus A$. Further, let $a = |A|$ and $b = |B|$, and so $n = a + b$. Since $T \not\cong T_k$ and $T \not\cong T^*_k$, we must have $k_1+k_2 \le k-1$. We note that $1 \le a \le 2k - 1$. Let $C_1 = C \cap A$. Let us first assume that $a = 1$. In this case, we have $k_1+k_2 = 0$ and so, the tree $T$ has no attachments of Type-1 or Type-2. Now, by the construction of the canonical set $C$, we have $|C| \le \frac{5}{6} \le \left( \frac{2\Delta - 1}{2\Delta} \right) n$ and hence, the result holds in this case. Let us therefore assume now that $a \ge 2$. Since neither $\bk = (1,0,1,0,0,0)$ nor $\bk = (1,0,0,0,1,0)$, by construction of the canonical set $C$, we have $|C_1| = a - 1 \le \left( \frac{2k - 2}{2k - 1} \right)a$. Further, let $C_2 = C \cap B$, and so $|C_2| \le \frac{4}{5}b$. Since $C$ is an IO-code of $T$ by Observation~\ref{obser1} and $\Delta \ge 3$, we infer that
\[
\ioc(T) \le |C| = |C_1| + |C_2| \le \left( \frac{2k - 2}{2k - 1} \right)a + \frac{4}{5}b < \left( \frac{2\Delta - 1}{2\Delta} \right) n. \hspace*{0.25cm}
\]
\end{proof}

\subsection{Theorem~\ref{thm:main1} for trees}
\label{S:key}

In this section, we prove the following Theorem~\ref{thm:main-tree} which we will need when proving our main result, namely Theorem~\ref{thm:main1}, for graphs that contain no $4$-cycles. First, we need a lemma.

\begin{lemma} \label{lem:T1 U T2}
Let $T$ be an open twin-free tree on $n$ vertices and let $T_1$ and $T_2$ be open twin-free trees obtained by removing an edge from $T$. Moreover, let $S_i$ be an IO-code of $T_i$ for each $i \in \{1, 2\}$ and let $S = S_1 \cup S_2$. Then we have the following.
\begin{enumerate}
\item[(a)] The set $S$ is an IO-code of $T$.
\item[(b)] If $|S_i| \le \left( \frac{2\Delta-1}{2\Delta} \right) |V(T_i)|$, then $|S| \le \left( \frac{2\Delta-1}{2\Delta} \right) n$.
\end{enumerate}
\end{lemma}

\begin{proof} Let the edge removed from the tree $T$ be $v_1v_2$. Moreover, by a possible renaming of the vertices, let $v_1 \in V(T_1)$ and $v_2 \in V(T_2)$. 

\noindent (a) Suppose, to the contrary, that there exists a pair $x,y$ of vertices of $T$ which are not identified by $S$. Then we must have, without loss of generality, $x = v_1$ and $y \in V(T_2)$ with $v_2y \in E(T_2)$. Since $S_1$ is a TD-set of $T_1$, we must have $v_1z \in E(T_1)$ for some vertex $z \in V(T_1) \cap S_1$. Now since $S$ does not identify the pair $x,y$, it implies that $yz \in E(T)$ as well which further implies that $zxv_2yz$ is a $4$-cycle in $T$ which contradicts the fact that the latter is a tree. Hence, $S$ is an IO-code of $T$.

\noindent (b) Since the sets $S_1$ and $S_2$ are disjoint, we therefore have
\[
|S| = |S_1|+|S_2| \le \left( \frac{2\Delta-1}{2\Delta} \right) \left( |V(T_1)|+|V(T_2)| \right) = \left( \frac{2\Delta-1}{2\Delta} \right) n.
\]
\end{proof}

\begin{theorem}
\label{thm:main-tree}
For $\Delta \ge 3$ a fixed integer, if $T \not \cong T_\Delta$ is an open twin-free tree of order~$n \ge 5$ that satisfies $\Delta(T) \le \Delta$, then
\[
\ioc(T) \le \left( \frac{2\Delta - 1}{2\Delta} \right) n.
\]
On the other hand, when $T \cong  T_{\Delta}$, then we have $\ioc(T) = \left( \frac{2\Delta}{2\Delta + 1} \right)n$.
\end{theorem}
\begin{proof} Let $\Delta \ge 3$ be a fixed integer. We proceed by induction on the order~$n \ge 5$ of a tree $T$ that is open twin-free and satisfies $\Delta(T) \le \Delta$. Since $T$ has order~$n \ge 5$ and is twin-free, we note that $\diam(T) \ge 4$ and $T$ has no strong support vertices, that is, every support has exactly one leaf neighbor. We proceed further with a series of claims.

\begin{claim}
\label{claim1}
If $\diam(T) = 4$, then one of the following holds. \\ [-20pt]
\begin{enumerate}
\item[{\rm (a)}] $T \cong T_{\Delta}$ and $\ioc(T) = \left( \frac{2\Delta}{2\Delta + 1} \right)n$.  \1
\item[{\rm (b)}] $T \not \cong T_{\Delta}$ and $\ioc(T) \le \left( \frac{2\Delta - 1}{2\Delta} \right) n$.
\end{enumerate}
\end{claim}
\proofclaim Suppose that $\diam(T) = 4$. Since $T$ has no strong support vertex, either $T \cong T_k$ for some $k$ where $2 \le k \le \Delta$ or $T \cong T_k^*$ for some $k$ where $3 \le k \le \Delta$.
If $T \cong T_{\Delta}$, then by Proposition~\ref{prop:prop1}(a), we have $\ioc(T) = \left( \frac{2\Delta}{2\Delta + 1} \right)n$, and so statement~(a) of the claim holds. If $T \cong T_k$ for some $k$ where $2 \le k \le \Delta - 1$, then by Proposition~\ref{prop:prop1}(a),
\[
\ioc(T) = \left( \frac{2k}{2k + 1} \right)n \le \left( \frac{2\Delta - 2}{2\Delta - 1} \right)n < \left( \frac{2\Delta - 1}{2\Delta} \right) n.
\]

If $T \cong T_k^*$ for some $k$ where $3 \le k \le \Delta$, then by Proposition~\ref{prop:prop1}(b),
\[
\ioc(T) = \left( \frac{2k-1}{2k} \right)n \le \left( \frac{2\Delta-1}{2\Delta} \right)n,
\]
and so statement~(b) of the claim holds.~\smallqed

\medskip
By Claim~\ref{claim1}, we may assume that $\diam(T) \ge 5$, for otherwise the desired result holds. In what follows, let $e = v_1v_2$ be an edge of $T$ and let $F_i$ be the component of $T - e$ that contains the vertex $v_i$ for $i \in [2]$. Further, let $F_i$ have order~$n_i$ where $i \in [2]$, and so $n = n_1 + n_2$. We note that $T$ is obtained from the trees $F_1$ and $F_2$ by adding back the edge~$e$.

\begin{claim}
\label{claim2}
If $F_1 \cong T_k$ and $v_1$ is the central vertex of degree~$k$ in $F_1$ for some $k$ where $2 \le k \le \Delta - 1$, then
\[
\ioc(T) \le \left( \frac{2\Delta - 1}{2\Delta} \right) n.
\]
\end{claim}
\proofclaim Suppose that $F_1 \cong T_k$ and $v_1$ is the central vertex of degree~$k$ in $F_1$ for some $k$ where $2 \le k \le \Delta - 1$. Then we have $F_1 \cong T(v_1;0,k,0,0,0,0) \in \cT$. Let $S_1$ denote the canonical set of the tree $F_1$. Thus, $S_1$ contains all vertices of $F_1$, except for exactly one leaf, $u_1$ say. In particular, we note that $S_1$ contains the vertex $v_1$ and all support vertices of $F_1$. We first consider the case when $n_2 \le 4$. Since $T \not \cong T_\Delta$, we have $n_2 \ne 2$. Moreover, for $n_2 \in \{1,3,4\}$, the tree $T$ is isomorphic to some tree in the special family $\cT$. Therefore, in the case that $n_2 \le 4$, we are done by Proposition~\ref{prop2_new}. Hence, we may assume that $n_2 \ge 5$. 

Let us now assume that $F_2$ is open twin-free. Suppose that $F_2 \cong T_{\Delta}$. In this case, either $v_2$ is a support vertex in $F_2$ (as illustrated in Figure~\ref{fig:claim2-end}(a)) or $v_2$ is a leaf in $F_2$ (as illustrated in Figure~\ref{fig:claim2-end}(b)). In both cases, let $S_1 = V(F_1) \setminus \{u_1\}$. If $v_2$ is a support vertex in $F_2$, then let $S_2$ consist of all vertices of $F_2$ except for exactly two leaves, one of which is the leaf neighbor of $v_2$ in $F_2$. If $v_2$ is a leaf in $F_2$, then let $S_2$ consist of all vertices of $F_2$ except for exactly one leaf which is different from the vertex $v_2$. Let $S = S_1 \cup S_2$. In the first case when $v_2$ is a support vertex in $F_2$, the resulting set $S$ is illustrated by the shaded vertices in Figure~\ref{fig:claim2-end}(a). In the second case when $v_2$ is a leaf in $F_2$, the resulting set $S$ is illustrated by the shaded vertices in Figure~\ref{fig:claim2-end}(b). In both cases, the set $S$ is an IO-code of $T$ by Lemma~\ref{lem:T1 U T2}(a). Thus, since $|S| = |S_1| + |S_2| \le 2k + 2\Delta$ and $n = (2k+1) + (2\Delta + 1)$, and since $k \le \Delta - 1$, we infer that
\[
\ioc(T) \le |S| \le \left( \frac{2k + 2\Delta}{2k + 2\Delta + 2} \right) n = \left( \frac{k + \Delta}{k + \Delta + 1} \right) n \le \left( \frac{2\Delta - 1}{2\Delta} \right) n.
\]

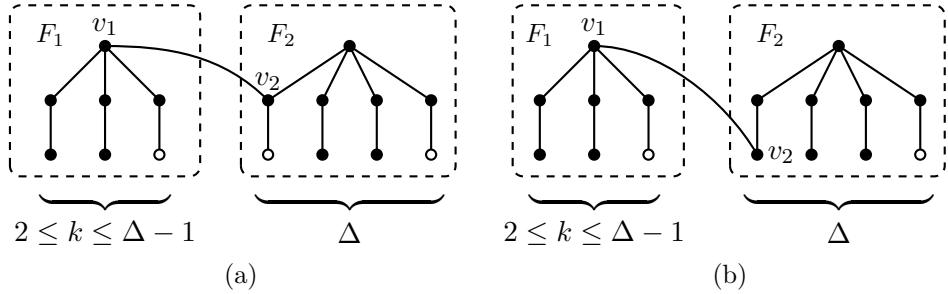
\begin{figure}[htb]
\begin{center}
\begin{tikzpicture}[scale=.85,style=thick,x=0.85cm,y=0.85cm]
\def\vr{2.25pt}
\path (0,0) coordinate (u1);
\path (0,1) coordinate (v1);
\path (1,0) coordinate (u2);
\path (1,1) coordinate (v2);
\path (2,0) coordinate (u3);
\path (2,1) coordinate (v3);
\path (1,2) coordinate (v);
\draw (v)--(v1)--(u1);
\draw (v)--(v2)--(u2);
\draw (v)--(v3)--(u3);
\path (4,0) coordinate (w1);
\path (4,1) coordinate (x1);
\path (5,0) coordinate (w2);
\path (5,1) coordinate (x2);
\path (6,0) coordinate (w3);
\path (6,1) coordinate (x3);
\path (7,0) coordinate (w4);
\path (7,1) coordinate (x4);
\path (5.5,2) coordinate (w);
\draw (w)--(x1)--(w1);
\draw (w)--(x2)--(w2);
\draw (w)--(x3)--(w3);
\draw (w)--(x4)--(w4);
\draw (v) to[out=0,in=135, distance=1cm] (x1);
\draw (1,-0.7) node {$\underbrace{\phantom{111111111}}$};
\draw (1,-1.45) node {$2 \le k \le \Delta - 1$};
\draw (5.5,-0.7) node {$\underbrace{\phantom{1111111111111}}$};
\draw (5.5,-1.45) node {$\Delta$};
\draw [style=dashed,rounded corners] (-0.75,-0.4) rectangle (2.75,2.75);
\draw (0,2.2) node {{\small $F_1$}};
\draw [style=dashed,rounded corners] (3.5,-0.4) rectangle (7.45,2.75);
\draw (4.25,2.2) node {{\small $F_2$}};
\draw (v1) [fill=black] circle (\vr);
\draw (v2) [fill=black] circle (\vr);
\draw (v3) [fill=black] circle (\vr);
\draw (u1) [fill=black] circle (\vr);
\draw (u2) [fill=black] circle (\vr);
\draw (u3) [fill=white] circle (\vr);
\draw (w1) [fill=white] circle (\vr);
\draw (w2) [fill=black] circle (\vr);
\draw (w3) [fill=black] circle (\vr);
\draw (w4) [fill=white] circle (\vr);
\draw (x1) [fill=black] circle (\vr);
\draw (x2) [fill=black] circle (\vr);
\draw (x3) [fill=black] circle (\vr);
\draw (x4) [fill=black] circle (\vr);
\draw (v) [fill=black] circle (\vr);
\draw (w) [fill=black] circle (\vr);
\draw[anchor = south] (v) node {$v_1$};
\draw[anchor = south] (x1) node {$v_2$};
\draw (3.5,-2.25) node {{\small (a)}};
\path (9,0) coordinate (u1);
\path (9,1) coordinate (v1);
\path (10,0) coordinate (u2);
\path (10,1) coordinate (v2);
\path (11,0) coordinate (u3);
\path (11,1) coordinate (v3);
\path (10,2) coordinate (v);
\draw (v)--(v1)--(u1);
\draw (v)--(v2)--(u2);
\draw (v)--(v3)--(u3);
\path (13,0) coordinate (w1);
\path (13,1) coordinate (x1);
\path (14,0) coordinate (w2);
\path (14,1) coordinate (x2);
\path (15,0) coordinate (w3);
\path (15,1) coordinate (x3);
\path (16,0) coordinate (w4);
\path (16,1) coordinate (x4);
\path (14.5,2) coordinate (w);
\draw (w)--(x1)--(w1);
\draw (w)--(x2)--(w2);
\draw (w)--(x3)--(w3);
\draw (w)--(x4)--(w4);
\draw (v) to[out=0,in=120, distance=1cm] (w1);
\draw (10,-0.7) node {$\underbrace{\phantom{111111111}}$};
\draw (10,-1.45) node {$2 \le k \le \Delta - 1$};
\draw (14.5,-0.7) node {$\underbrace{\phantom{1111111111111}}$};
\draw (14.5,-1.45) node {$\Delta$};
\draw [style=dashed,rounded corners] (8.5,-0.4) rectangle (11.65,2.75);
\draw (9,2.2) node {{\small $F_1$}};
\draw [style=dashed,rounded corners] (12.5,-0.4) rectangle (16.45,2.75);
\draw (13.25,2.2) node {{\small $F_2$}};
\draw (v1) [fill=black] circle (\vr);
\draw (v2) [fill=black] circle (\vr);
\draw (v3) [fill=black] circle (\vr);
\draw (u1) [fill=black] circle (\vr);
\draw (u2) [fill=black] circle (\vr);
\draw (u3) [fill=white] circle (\vr);
\draw (w1) [fill=black] circle (\vr);
\draw (w2) [fill=black] circle (\vr);
\draw (w3) [fill=black] circle (\vr);
\draw (w4) [fill=white] circle (\vr);
\draw (x1) [fill=black] circle (\vr);
\draw (x2) [fill=black] circle (\vr);
\draw (x3) [fill=black] circle (\vr);
\draw (x4) [fill=black] circle (\vr);
\draw (v) [fill=black] circle (\vr);
\draw (w) [fill=black] circle (\vr);
\draw[anchor = south] (v) node {$v_1$};
\draw[anchor = west] (w1) node {$v_2$};
\draw (12.5,-2.25) node {{\small (b)}};
\end{tikzpicture}
\caption{Possible trees in the proof of Claim~\ref{claim2}.}
\label{fig:claim2-end}
\end{center}
\end{figure}

Hence we may assume that $F_2 \not \cong T_{\Delta}$. Applying the inductive hypothesis to the open twin-free tree $F_2$ of order~$n_2 \ge 5$, there exists an IO-code, $S_2$ say, of $F_2$ such that $\ioc(F_2) = |S_2| \le \left( \frac{2\Delta - 1}{2\Delta} \right) n_2$. Let $S = S_1 \cup S_2$. Since the vertex $v_2$ in $F_2$ is totally dominated by some other vertex of $F_2$ in $S_2$, the set $S$ is an IO-code of the tree $T$. By Proposition~\ref{prop2_new}, we have $|S_1| \le \left(\frac{2\Delta - 1}{2\Delta}\right) n_1$, and so
\[
\ioc(T) \le |S| = |S_1| + |S_2| \le \left( \frac{2\Delta - 1}{2\Delta} \right) n_1 + \left( \frac{2\Delta - 1}{2\Delta} \right) n_2 = \left( \frac{2\Delta - 1}{2\Delta} \right)n.
\]

Let us now assume that $F_2$ contains open twins. Since the tree $T$ contains no open twins, we infer that $v_2$ is a leaf in $F_2$ and is an open twin in $F_2$ with some other leaf, say $u_2$, in $F_2$. Let $v$ be the common neighbor of $v_2$ and $u_2$ in $F_2$. Let $T' = F_2 - v_2$, and let $T'$ have order~$n'$, and so $n' = n_2 - 1 \ge 4$. Since $T$ is open twin-free, so too is the tree $T'$.

Suppose that $n' = 4$, and so $n_2 = 5$. In this case, the tree $T \cong T(v_1;0,k,0,0,0,1) \in \cT$ and therefore, we are done by Proposition~\ref{prop2_new}. Hence we may assume that $n' \ge 5$. Suppose that $T' \cong T_{\Delta}$. In this case, the vertex $v$ is a support vertex in the tree $T'$ and the tree $T$ is illustrated in Figure~\ref{fig:claim2-middle}. Let $u_2$ denote the leaf neighbor of the vertex~$v$. We note that $n_1 = 2k+1$ and $n_2 = 2\Delta + 2$, and so $n = 2k + 2\Delta + 3$. We now let $S_2$ consist of all vertices in $T'$, except for the leaf $u_2$. Thus, $|S_2| = 2\Delta$. Further, we let $S = S_1 \cup S_2$, and so the set $S$ is indicated by the shaded vertices in Figure~\ref{fig:claim2-middle}. Then by Lemma~\ref{lem:T1 U T2}, the set $S$ is an IO-code of the tree, and so
\[
\ioc(T) \le |S| = 2k + 2\Delta
= \left( \frac{2k + 2\Delta}{2k + 2\Delta + 3} \right) n
\le \left( \frac{4\Delta-2}{4\Delta + 1} \right) n
< \left( \frac{2\Delta - 1}{2\Delta} \right) n.
\]

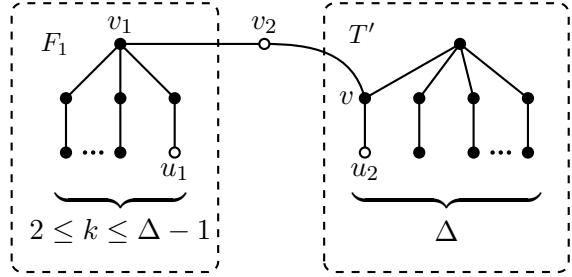
\begin{figure}[htb]
\begin{center}
\begin{tikzpicture}[scale=.85,style=thick,x=0.85cm,y=0.85cm]
\def\vr{2.25pt}
\def\vrs{0.5pt}
\path (0,0) coordinate (u1);
\path (0,1) coordinate (v1);
\path (1,0) coordinate (u2);
\path (1,1) coordinate (v2);
\path (2,0) coordinate (u3);
\path (2,1) coordinate (v3);
\path (3.65,2) coordinate (z);
\path (1,2) coordinate (v);
\path (0.35,0) coordinate (p1);
\path (0.5,0) coordinate (p2);
\path (0.65,0) coordinate (p3);
\draw (v)--(v1)--(u1);
\draw (v)--(v2)--(u2);
\draw (v)--(v3)--(u3);
\draw (v)--(z);
\path (5.5,0) coordinate (w1);
\path (5.5,1) coordinate (x1);
\path (6.5,0) coordinate (w2);
\path (6.5,1) coordinate (x2);
\path (7.5,0) coordinate (w3);
\path (7.5,1) coordinate (x3);
\path (8.5,0) coordinate (w4);
\path (8.5,1) coordinate (x4);
\path (7.25,2) coordinate (w);
\path (7.85,0) coordinate (q1);
\path (8,0) coordinate (q2);
\path (8.15,0) coordinate (q3);

\draw (w)--(x1)--(w1);
\draw (w)--(x2)--(w2);
\draw (w)--(x3)--(w3);
\draw (w)--(x4)--(w4);
\draw (z) to[out=0,in=105, distance=0.75cm] (x1);
\draw (1,-0.7) node {$\underbrace{\phantom{111111111}}$};
\draw (1,-1.45) node {$2 \le k \le \Delta - 1$};
\draw (7,-0.7) node {$\underbrace{\phantom{1111111111111}}$};
\draw (7,-1.45) node {$\Delta$};
\draw [style=dashed,rounded corners] (-1,-2.2) rectangle (2.8,2.75);
\draw (-0.2,2) node {{\small $F_1$}};
\draw [style=dashed,rounded corners] (4.75,-2.2) rectangle (9.25,2.75);
\draw (5.45,2.2) node {{\small $T'$}};
\draw (v1) [fill=black] circle (\vr);
\draw (v2) [fill=black] circle (\vr);
\draw (v3) [fill=black] circle (\vr);
\draw (u1) [fill=black] circle (\vr);
\draw (u2) [fill=black] circle (\vr);
\draw (u3) [fill=white] circle (\vr);
\draw (w1) [fill=white] circle (\vr);
\draw (w2) [fill=black] circle (\vr);
\draw (w3) [fill=black] circle (\vr);
\draw (w4) [fill=black] circle (\vr);
\draw (x1) [fill=black] circle (\vr);
\draw (x2) [fill=black] circle (\vr);
\draw (x3) [fill=black] circle (\vr);
\draw (x4) [fill=black] circle (\vr);
\draw (v) [fill=black] circle (\vr);
\draw (w) [fill=black] circle (\vr);
\draw (z) [fill=white] circle (\vr);
\draw (p1) [fill=black] circle (\vrs);
\draw (p2) [fill=black] circle (\vrs);
\draw (p3) [fill=black] circle (\vrs);
\draw (q1) [fill=black] circle (\vrs);
\draw (q2) [fill=black] circle (\vrs);
\draw (q3) [fill=black] circle (\vrs);
\draw[anchor = south] (v) node {$v_1$};
\draw[anchor = south] (z) node {$v_2$};
\draw[anchor = east] (x1) node {$v$};
\draw[anchor = north] (w1) node {$u_2$};
\draw[anchor = north] (u3) node {$u_1$};
\end{tikzpicture}
\caption{A possible tree in the proof of Claim~\ref{claim2}.}
\label{fig:claim2-middle}
\end{center}
\end{figure}

Hence, we may assume that $T' \not \cong T_{\Delta}$. Recall that $n' \ge 5$. Applying the inductive hypothesis to the open twin-free tree $T'$ of order~$n' \ge 5$, there exists an IO-code, $S'$ say, of $T'$ such that $\ioc(T') = |S'| \le \left( \frac{2\Delta - 1}{2\Delta} \right) n'$. Now let $S = S_1 \cup S'$. Since $u_2$ is a leaf in $T'$, the IO-code $S'$ of $T'$ contains the support vertex $v$, which has two neighbors in $S$. Therefore, the set $S$ is an IO-code of the tree $T$, and so
\[
\ioc(T) \le |S| = |S_1| + |S'| \le \left( \frac{2\Delta - 1}{2\Delta} \right) n_1 + \left( \frac{2\Delta - 1}{2\Delta} \right) n' < \left( \frac{2\Delta - 1}{2\Delta} \right)n.
\]
This completes the proof of Claim~\ref{claim2}.~\smallqed

\medskip
By Claim~\ref{claim2}, we may assume that the removal of any edge from $T$ does not produce a component $T'$ isomorphic to $T_k$, where $2 \le k \le \Delta - 1$ and where the central vertex of the subdivided star $T'$ is incident with the deleted edge, for otherwise the desired result follows. As a consequence of Claim~\ref{claim2}, we have the following property of the tree $T$.

\begin{claim}
\label{claim2b}
If $T$ contains an edge whose removal produces a component isomorphic to $T_{\Delta}$, then
\[
\ioc(T) \le \left( \frac{2\Delta - 1}{2\Delta} \right) n.
\]
\end{claim}
\proofclaim Suppose that $T$ contains an edge $f = xy$ whose removal produces a component, $T'$ say, isomorphic to $T_{\Delta}$. Renaming $x$ and $y$ if necessary, we may assume that $x \in V(T')$. Let $v$ be the central vertex of the subdivided star in $T'$, and so $v$ has degree~$\Delta$ in $T'$ (and in the original tree $T$). If $x$ is a support vertex of $T'$, then we infer from Claim~\ref{claim2} with $v_1 = v$ and $v_2 = x$ that $\ioc(T) \le \left( \frac{2\Delta - 1}{2\Delta} \right) n$. If $x$ is a leaf of $T'$, then we let $w$ denote the common neighbor of $v$ and $x$, and we infer from Claim~\ref{claim2} with $v_1 = v$ and $v_2 = w$ that $\ioc(T) \le \left( \frac{2\Delta - 1}{2\Delta} \right) n$.~\smallqed

\medskip
By Claim~\ref{claim2b}, we may assume that the removal of any edge from $T$ does not produce a component isomorphic to $T_{\Delta}$, for otherwise the desired result follows.

For the rest of the proof, let $v_0v_1v_2 \ldots v_d$ to be a fixed longest path in $T$. Then $v_0$ and $v_d$ are leaves in $T$. Let us also assume from here on that the tree $T$ is rooted at the leaf $v_d$. Since $diam(T) \ge 5$ by assumption, we must have $d \ge 5$. This implies that,  for all $i \in [5]$, the vertex $v_i$ is the ancestor of $v_{i-1}$ in $T$. If $v$ is a vertex in $T$, then we denote by $T_v$ the maximal subtree rooted at $v$, and so $T_v$ consists of $v$ and all descendants of the vertex~$v$. Since $T$ has no strong support vertex, we note that $\deg_T(v_1) = 2$.

\begin{claim} \label{claim_T_vi}
If $\deg_T(v_2) \ge 4$, or if $\deg_T(v_i) \ge 3$ for $i \in \{3,4\}$, then
\[
\ioc(T) \le \left( \frac{2\Delta - 1}{2\Delta} \right) n.
\]
\end{claim}
\proofclaim
We first assume that $\deg_T(v_2) \ge 4$. Then, the maximal subtree $T_{v_2}$ of $T$ rooted at $v_2$ is isomorphic to a tree in the special family $\cT$. By Claim~\ref{claim2}, we may assume that $T_{v_2} \not \cong T_k$ for any $k \le \Delta - 1$. Hence, $T_{v_2}$ must be isomorphic to the reduced subdivided tree $T^*_k$. Let $T_{v_2}$ be on $n_1$ vertices and let $C$ be the canonical set of $T_{v_2}$. Then, by Observation~\ref{obser1}, the set $C$ is an IO-code of the tree $T_{v_2}$ and by Proposition~\ref{prop:prop1}, we have $|C| \le \left(\frac{2\Delta - 1}{2\Delta}\right)n_1$. Moreover, let $T'$ be the component of $T - v_2v_3$ that contains the vertex $v_3$, and let $T'$ have order~$n'$. If $n' \le 4$, then $T$ is isomorphic to a tree in the special family $\cT$ and thus, we are done by Proposition~\ref{prop2_new}. Let us therefore assume that $n' \ge 5$. By Claim~\ref{claim2b}, we also assume that $T'$ is not isomorphic to $T_\Delta$.

First, we assume that $T'$ is open twin-free. Therefore, applying the inductive hypothesis to the open twin-free tree $T'$, there exists an IO-code, $S'$ say, of $T'$ such that $\ioc(T') = |S'| \le \left( \frac{2\Delta - 1}{2\Delta} \right) n'$. Let $S = C \cup S'$. Then by Lemma~\ref{lem:T1 U T2}, the set $S$ is an IO-code of the tree $T$ with 
\[ 
\ioc(T) \le |S| \le \left( \frac{2\Delta - 1}{2\Delta} \right)n.
\]

We now suppose that $T'$ contains open twins. Since the tree $T$ contains no open twins, we infer that $v_3$ is a leaf in $T'$ and is an open twin in $T'$ with some other leaf, say $u_3$, in $T'$. We note that $v_4$ is the common neighbor of $v_3$ and $u_3$. Let $T'' = T' - v_3$, and let $T''$ have order~$n''$, and so $n'' = n' - 1 \ge 4$. If $n'' = 4$, since $v_3$ and $u_3$ are open twins in $T''$, it implies that $T \cong T(v_2;1,k-2,0,0,0,1) \in \cT$ for some $k = \deg_T(v_2)-1 \le \Delta - 1$. Hence, in this case, we are again done by Proposition~\ref{prop2_new}. Let us therefore assume that $n'' \ge 5$. Since $T$ is open twin-free, so too is the tree $T''$. By Claim~\ref{claim2b}, we infer that $T'' \not \cong T_{\Delta}$. Applying the inductive hypothesis to the open twin-free tree $T''$ of order~$n'' \ge 5$, there exists an IO-code, $S''$ say, of $T''$ such that $\ioc(T'') = |S''| \le \left( \frac{2\Delta - 1}{2\Delta} \right) n''$. Let $S = C \cup S''$. We note that since $u_3$ is a leaf, the vertex $v_4 \in S''$ and so the vertex $v_3$ has two neighbors in the set $S$ and thus is identified by $S$. Hence, the set $S = C \cup S''$ is an IO-code of $T$, and so
\[
\begin{array}{lcl}
\ioc(T) \le |S| & = & |C| + |S''| \2 \\
& \le & \left( \frac{2\Delta - 1}{2\Delta} \right) n_1 + \left( \frac{2\Delta - 1}{2\Delta} \right) n'' \2 \\
& < & \left( \frac{2\Delta - 1}{2\Delta} \right)(n_1 + 1) + \left( \frac{2\Delta - 1}{2\Delta} \right) n'' \2 \\
& = & \left( \frac{2\Delta - 1}{2\Delta} \right)(n - n'') + \left( \frac{2\Delta - 1}{2\Delta} \right) n'' \2 \\
& = & \left( \frac{2\Delta - 1}{2\Delta} \right)n.
\end{array}
\]

This proves the result in the case that $\deg_T(v_2) \ge 4$. We can therefore assume from now on that $\deg_T(v_2) \le 3$, or else, the result holds. In fact, more strongly, we can also assume that whenever $\deg_T(v_2) = 3$, the subtree $T_{v_2}$ is not isomorphic to the subdivided star $T_2$. This is due to our previous assumption using Claim~\ref{claim2} that $T_{v_2} \not \cong T_k$ for any $k \le \Delta - 1$.

We now look at the case that $\deg_T(v_i) \ge 3$ for all $i \in \{3,4\}$. As before, let $n_1$ be the number of vertices of $T_{v_i}$, the subtree of $T$ rooted at $v_i$. Then $T_{v_3}$ (respectively, $T_{v_4}$) is isomorphic to a tree $T(r;\bk)$ in $\cT$ such that $k_3+k_5 \ge 1$ (respectively, $k_4+k_6 \ge 1$). This implies that $T_{v_i}$ is neither a subdivided star nor a reduced subdivided star. Let $C$ be the canonical set of $T_{v_i}$. Then by Observation~\ref{obser1}, the set $C$ is an IO-code of the tree $T_{v_i}$. Moreover, by Proposition~\ref{prop2_new}, we have $|C| \le \left(\frac{2\Delta - 1}{2\Delta}\right)n_1$. Now, taking $T'$ to be the component of $T - v_3v_4$ (respectively, $T - v_4v_5$) containing the vertex $v_4$ (respectively, $v_5$), the result follows in this case by the same analysis on the tree $T'$ as in the case $\deg_T(v_2) \ge 4$.~\smallqed

Hence, by Claim~\ref{claim_T_vi}, we can assume from here on that $2 \le \deg_T(v_2) \le 3$ and $\deg_T(v_i) = 2$ for all $i \in \{3,4\}$, or else, the desired result is achieved. We now look at the maximal subtree $T_{v_4}$ of $T$ rooted at $v_4$. By our earlier assumptions, either $T_{v_4}$ is a path $P_5$ with $v_4$ as one of the ends of the path as illustrated in Figure~\ref{fig:final}(a), or $T_{v_4}$ is a reduced subdivided star $T_{3}^*$ as illustrated in Figure~\ref{fig:final}(b) with $v_2$ as the central vertex of the reduced subdivided star. In the latter case, we let $u_1$ be the leaf neighbor of $v_2$.

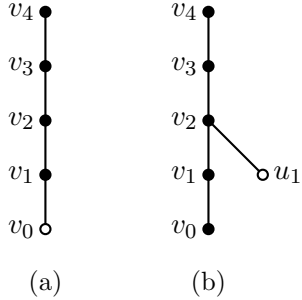
\begin{figure}[htb]
\begin{center}
\begin{tikzpicture}[scale=.85,style=thick,x=0.85cm,y=0.85cm]
\def\vr{2.25pt}
\path (2,0) coordinate (v0);
\path (2,1) coordinate (v1);
\path (2,2) coordinate (v2);
\path (2,3) coordinate (v3);
\path (2,4) coordinate (v4);
\draw (v4)--(v3)--(v2)--(v1)--(v0);
%
\draw (v1) [fill=black] circle (\vr);
\draw (v2) [fill=black] circle (\vr);
\draw (v3) [fill=black] circle (\vr);
\draw (v0) [fill=white] circle (\vr);
\draw (v4) [fill=black] circle (\vr);
\draw[anchor = east] (v4) node {$v_4$};
\draw[anchor = east] (v0) node {$v_0$};
\draw[anchor = east] (v1) node {$v_1$};
\draw[anchor = east] (v2) node {$v_2$};
\draw[anchor = east] (v3) node {$v_3$};
\draw (2,-1) node {{\small (a)}};
\path (5,0) coordinate (v0);
\path (5,1) coordinate (v1);
\path (5,2) coordinate (v2);
\path (6,1) coordinate (u1);
\path (5,3) coordinate (v3);
\path (5,4) coordinate (v4);
\draw (v4)--(v3)--(v2)--(v1)--(v0);
\draw (v2)--(u1);
%
\draw (v1) [fill=black] circle (\vr);
\draw (v2) [fill=black] circle (\vr);
\draw (v3) [fill=black] circle (\vr);
\draw (v0) [fill=black] circle (\vr);
\draw (v4) [fill=black] circle (\vr);
\draw (u1) [fill=white] circle (\vr);
\draw[anchor = east] (v4) node {$v_4$};
\draw[anchor = east] (v0) node {$v_0$};
\draw[anchor = east] (v1) node {$v_1$};
\draw[anchor = west] (u1) node {$u_1$};
\draw[anchor = east] (v2) node {$v_2$};
\draw[anchor = east] (v3) node {$v_3$};
\draw (5,-1) node {{\small (b)}};
\end{tikzpicture}
\caption{Two possible subtrees in $T$.}
\label{fig:final}
\end{center}
\end{figure}

If $T_{v_4}$ is a path $P_5$ with $v_4$ as one of the ends of the path, then we let $S_1 = \{v_1,v_2,v_3,v_4\}$. In this case, the set $S_1$, indicated by the shaded vertices in Figure~\ref{fig:final}(a), is an IO-code of $T_{v_4}$. If $T_{v_4}$ is a reduced subdivided star $T_{3}^*$, then we let $S_1 = \{v_0,v_1,v_2,v_3,v_4\}$. In this case, the set $S_1$, indicated by the shaded vertices in Figure~\ref{fig:final}(b), is an IO-code of $T_{v_4}$.

Let $T'$ be the component of $T - v_4v_5$ that contains the vertex $v_5$, and let $T'$ have order~$n'$. If $n' \le 4$, then the tree $T$ is determined and is isomorphic to a tree in $\cT$, and hence, we infer the desired upper bound by Proposition~\ref{prop2_new}. We may therefore assume that $n' \ge 5$.

Suppose now that $T'$ is open twin-free. By our earlier assumptions, we infer that $T' \not \cong T_{\Delta}$. Applying the inductive hypothesis to the open twin-free tree $T'$, there exists an IO-code, $S'$ say, of $T'$ such that $\ioc(T') = |S'| \le \left( \frac{2\Delta - 1}{2\Delta} \right) n'$. Let $S = S_1 \cup S'$. Then by Lemma~\ref{lem:T1 U T2}(a), the set $S$ is an IO-code of the tree $T$. If $T'$ is a path $P_5$, then $n_1 = 5$ and $|S_1| = 4 = \frac{4}{5}n_1 = \frac{4}{5}(n - n')$. If $T'$ is a reduced subdivided star $T_{3}^*$, then $n_1 = 6$ and $|S_1| = 5 = \frac{5}{6}n_1 = \frac{5}{6}(n - n')$. In both cases, $|S_1| \le \frac{5}{6}(n - n') \le \left( \frac{2\Delta - 1}{2\Delta} \right) (n-n')$. Therefore, using Lemma~\ref{lem:T1 U T2}(b), we have
\[
\ioc(T) \le |S| \le \left( \frac{2\Delta - 1}{2\Delta} \right)n.
\]

Suppose now that $T'$ contains open twins. Since the tree $T$ contains no open twins, we infer that $v_5$ is a leaf in $T'$ and is an open twin in $T'$ with some other leaf, say $u_5$, in $T'$. We note that $v_6$ is the common neighbor of $v_5$ and $u_5$. Let $T'' = T' - v_5$, and let $T''$ have order~$n''$, and so $n'' = n' - 1 \ge 4$. If $n'' = 4$, then $T$ is isomorphic to $T(v_4;0,0,0,1,0,1)$ if $T_{v_4} \cong P_5$ and to $T(v_4;0,0,0,0,0,2)$ if $T_{v_4} \cong T^*_3$. In both cases, the result follows by Proposition~\ref{prop2_new}. Therefore, let us assume that $n'' \ge 5$. Since $T$ is open twin-free, so too is the tree $T''$. By Claim~\ref{claim2b}, we infer that $T'' \not \cong T_{\Delta}$. Applying the inductive hypothesis to the open twin-free tree $T''$ of order~$n'' \ge 5$, there exists an IO-code, $S''$ say, of $T''$ such that $\ioc(T'') = |S''| \le \left( \frac{2\Delta - 1}{2\Delta} \right) n''$. Let $S = S_1 \cup S''$. We note that since $u_5$ is a leaf, the vertex $v_6 \in S''$, and so $v_5$ is identified by the set $S$. Therefore, the set $S$ is an IO-code of the tree $T$. If $T_{v_4}$ is a path $P_5$, then $n_1 = 5$ and $|S_1| = 4 = \frac{2}{3}(n_1 + 1) = \frac{2}{3}(n - n'')$. If $T_{v_4}$ is a reduced subdivided star $T_{3}^*$, then $n_1 = 6$ and $|S_1| = 5 = \frac{5}{7}(n_1 + 1) = \frac{5}{7}(n - n'')$. In both cases, $|S_1| \le \frac{5}{7}(n - n'')$. Therefore,
\[
\begin{array}{lcl}
\ioc(T) \le |S| & = & |S_1| + |S''| \1 \\
& \le & \frac{5}{7} (n - n'') + \left( \frac{2\Delta - 1}{2\Delta} \right) n'' \2 \\
& < & \left( \frac{2\Delta - 1}{2\Delta} \right)(n - n'') + \left( \frac{2\Delta - 1}{2\Delta} \right) n'' \2 \\
& = & \left( \frac{2\Delta - 1}{2\Delta} \right)n.
\end{array}
\]
This completes the proof of Theorem~\ref{thm:main-tree}.
\end{proof}

\section{Proof of Theorem~\ref{thm:main1}}
\label{S:proof}

In this section, we present a proof of our main result, namely Theorem~\ref{thm:main1}. We first present the following preliminary result.

\begin{proposition}
\label{prop:main2}
For $\Delta \ge 3$ a fixed integer, if $G$ is a graph of order~$n$ that is isomorphic to a graph obtained from a subdivided star $T_k$, where $2 \le k \le \Delta$, by adding an edge $e$ to $T_k$ in such a way that $G \cong T_k + e$ contains no $4$-cycle and satisfies $\Delta(G) \le \Delta$, then
\[
\ioc(G) \le \left( \frac{2\Delta - 1}{2\Delta} \right) n.
\]
\end{proposition}
\begin{proof}
It can be verified that the graph $G$ in the statement of the proposition must be isomorphic to exactly one of the three graphs $G_1$, $G_2$ and $G_3$ portrayed in Figure~\ref{fig:prop:main2}. In each of these cases, we have $n = 2k+1$.

\begin{figure}[htb]
\begin{center}
\begin{tikzpicture}[scale=.85,style=thick,x=0.85cm,y=0.85cm]
\def\vr{2.25pt}
\path (0,0) coordinate (u1);
\path (0,1) coordinate (v1);
\path (1,0) coordinate (u2);
\path (1,1) coordinate (v2);
\path (2,0) coordinate (u3);
\path (2,1) coordinate (v3);
\path (3,0) coordinate (u4);
\path (3,1) coordinate (v4);
\path (5,0) coordinate (u5);
\path (5,1) coordinate (v5);
\path (2.5,2) coordinate (v);
%
\draw (v)--(v1)--(u1);
\draw (v)--(v2)--(u2);
\draw (v)--(v3)--(u3);
\draw (v)--(v4)--(u4);
\draw (v)--(v5)--(u5);
\draw (v1)--(v2);
\draw (v1) [fill=black] circle (\vr);
\draw (v2) [fill=black] circle (\vr);
\draw (v3) [fill=black] circle (\vr);
\draw (v4) [fill=black] circle (\vr);
\draw (v5) [fill=black] circle (\vr);
\draw (u1) [fill=white] circle (\vr);
\draw (u2) [fill=white] circle (\vr);
\draw (u3) [fill=black] circle (\vr);
\draw (u4) [fill=black] circle (\vr);
\draw (u5) [fill=black] circle (\vr);
\draw (v) [fill=black] circle (\vr);
\draw (2.5,-2.5) node {{\small (a) $G_1$}};
\draw (4,1) node {$\cdots$};
\draw (4,0) node {$\cdots$};
\draw (2.5,-0.5) node {$\underbrace{\phantom{1111111111111111111}}$};
\draw (2.5,-1.2) node {$2 \le k \le \Delta$};
\draw (0.5,0.7) node {$e$};
\path (7,0) coordinate (u1);
\path (7,1) coordinate (v1);
\path (8,0) coordinate (u2);
\path (8,1) coordinate (v2);
\path (9,0) coordinate (u3);
\path (9,1) coordinate (v3);
\path (10,0) coordinate (u4);
\path (10,1) coordinate (v4);
\path (12,0) coordinate (u5);
\path (12,1) coordinate (v5);
\path (9.5,2) coordinate (v);
%
\draw (v)--(v1)--(u1);
\draw (v)--(v2)--(u2);
\draw (v)--(v3)--(u3);
\draw (v)--(v4)--(u4);
\draw (v)--(v5)--(u5);
\draw (u1)--(u2);
\draw (v1) [fill=white] circle (\vr);
\draw (v2) [fill=black] circle (\vr);
\draw (v3) [fill=black] circle (\vr);
\draw (v4) [fill=black] circle (\vr);
\draw (v5) [fill=black] circle (\vr);
\draw (u1) [fill=white] circle (\vr);
\draw (u2) [fill=black] circle (\vr);
\draw (u3) [fill=black] circle (\vr);
\draw (u4) [fill=black] circle (\vr);
\draw (u5) [fill=black] circle (\vr);
\draw (v) [fill=black] circle (\vr);
\draw (9.5,-2.5) node {{\small (b) $G_2$}};
\draw (11,1) node {$\cdots$};
\draw (11,0) node {$\cdots$};
\draw (9.5,-0.5) node {$\underbrace{\phantom{1111111111111111111}}$};
\draw (9.5,-1.2) node {$2 \le k \le \Delta$};
\draw (7.5,0.3) node {$e$};
\path (14,0) coordinate (u1);
\path (14,1) coordinate (v1);
\path (15,0) coordinate (u2);
\path (15,1) coordinate (v2);
\path (16,0) coordinate (u3);
\path (16,1) coordinate (v3);
\path (17,0) coordinate (u4);
\path (17,1) coordinate (v4);
\path (19,0) coordinate (u5);
\path (19,1) coordinate (v5);
\path (16.5,2) coordinate (v);
%
\draw (v)--(v1)--(u1);
\draw (v)--(v2)--(u2);
\draw (v)--(v3)--(u3);
\draw (v)--(v4)--(u4);
\draw (v)--(v5)--(u5);
\draw (v) to[out=0,in=40, distance=2cm] (u5);
\draw (v1) [fill=black] circle (\vr);
\draw (v2) [fill=black] circle (\vr);
\draw (v3) [fill=black] circle (\vr);
\draw (v4) [fill=black] circle (\vr);
\draw (v5) [fill=black] circle (\vr);
\draw (u1) [fill=white] circle (\vr);
\draw (u2) [fill=black] circle (\vr);
\draw (u3) [fill=black] circle (\vr);
\draw (u4) [fill=black] circle (\vr);
\draw (u5) [fill=black] circle (\vr);
\draw (v) [fill=black] circle (\vr);
\draw (16.5,-2.5) node {{\small (c) $G_3$}};
\draw (18,1) node {$\cdots$};
\draw (18,0) node {$\cdots$};
\draw (16.5,-0.5) node {$\underbrace{\phantom{1111111111111111111}}$};
\draw (16.5,-1.2) node {$2 \le k \le \Delta - 1$};
\draw (19.9,1.5) node {$e$};
\end{tikzpicture}
\caption{Possible graphs $G \cong T_k + e$ as in Proposition~\ref{prop:main2}.}
\label{fig:prop:main2}
\end{center}
\end{figure}
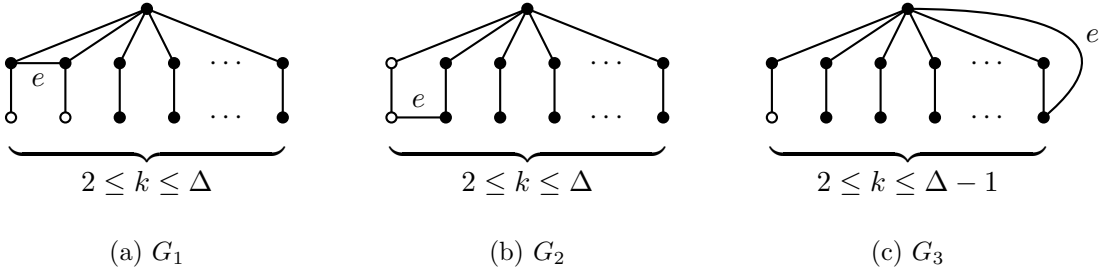

Suppose that $G \cong G_1$. In this case, $\ioc(G) \le 2k-1$, where the shaded vertices in Figure~\ref{fig:prop:main2}(a) form an IO-code of $G$. Thus, $\ioc(G) \le \left( \frac{2k - 1}{2k + 1} \right)n < \left( \frac{2\Delta - 1}{2\Delta} \right)n$.

Suppose that $G \cong G_2$. If $k = 2$, then $G \cong C_5$, $n = 5$ and $\ioc(G) = 4 = \frac{4}{5}n < \left( \frac{2\Delta - 1}{2\Delta} \right)n$. If $3 \le k \le \Delta$, then $\ioc(G) \le 2k-1$, where the shaded vertices in Figure~\ref{fig:prop:main2}(b) form an IO-code of $G$. Thus, $\ioc(G) \le \left( \frac{2k - 1}{2k + 1} \right)n < \left( \frac{2\Delta - 1}{2\Delta} \right)n$.

Suppose that $G \cong G_3$. In this case, $k \le \Delta - 1$ noting that $\Delta(G) \le \Delta$. Therefore we have $\ioc(G) \le 2k-1$, where the shaded vertices in Figure~\ref{fig:prop:main2}(c) form an IO-code of $G$. Thus, $\ioc(G) \le \left( \frac{2k - 1}{2k + 1} \right)n < \left( \frac{2\Delta - 1}{2\Delta} \right)n$.
\end{proof}

\medskip
We are now in a position to present a proof of our main result. Recall its statement.

\medskip
\noindent \textbf{Theorem~\ref{thm:main1}}. \emph{For $\Delta \ge 3$ a fixed integer, if $G \not \cong T_\Delta$ is an open twin-free connected graph of order~$n \ge 5$ that contains no $4$-cycles and satisfies $\Delta(G) \le \Delta$, then
\[
\ioc(G) \le \left( \frac{2\Delta - 1}{2\Delta} \right) n,
\]
On the other hand, when $G \cong  T_{\Delta}$, then we have $\ioc(G) = \left( \frac{2\Delta}{2\Delta + 1} \right)n$.
}

\begin{proof} We proceed by induction on the size~$m \ge 4$ of the connected graph $G$. If $m = 4$, then either $G \cong P_5$ or $G$ is isomorphic to a \emph{paw}, that is a complete graph on three vertices (equivalently, a $K_3$) with a leaf adjacent to one of the vertices of the $K_3$. Then, $\ioc(G) = 4$ if $G \cong P_5$ and $\ioc(G) = 3$ if $G$ is isomorphic to a paw. In both cases, we have $\ioc(G) \le \left( \frac{2\Delta-1}{2\Delta} \right)n$. This establishes the base case. For the inductive hypothesis, let $m \ge 5$ and assume that if $G'$ is a connected graph $G'$ of order~$n' \ge 5$ and size less than~$m$ that contains no $4$-cycle and is open twin-free and satisfies $\Delta(G') \le \Delta$, then $\ioc(G') \le \left( \frac{2\Delta - 1}{2\Delta} \right) n'$, unless $G' \cong  T_{\Delta}$.

We now consider the connected graph $G$ of order~$n \ge 5$ and size~$m$ that contains no $4$-cycles and is open twin-free and satisfies $\Delta(G) \le \Delta$. If the graph $G$ is a tree, then the desired result follows from Theorem~\ref{thm:main-tree}. Hence we may assume that the connected graph $G$ contains a cycle. If $m = 5$, then the graph $G$ is isomorphic to one of the graphs $G_1$, $G_2$ and $G_3$ illustrated in Figure~\ref{fig:prop:main2} and with $k=2$, and the desired result holds by Proposition~\ref{prop:main2}. Hence we may assume that $m \ge 6$.

Let $C$ be an induced cycle in $G$ and let $e = uv$ be an edge of the cycle $C$. Let $G' = G - e$ be the connected graph obtained from $G$ by deleting the cycle edge $e$. The resulting graph $G'$ has order~$n \ge 5$ and maximum degree $\Delta(G') \le \Delta(G) \le \Delta$. Since $G$ contains no $4$-cycles, neither does the graph $G'$. If $G' \cong T_\Delta$, then the graph $G$ satisfies the statement of Proposition~\ref{prop:main2}. Hence we may assume that $G' \not\cong T_\Delta$, for otherwise the desired result follows by Proposition~\ref{prop:main2}.

Suppose that $G'$ is open twin-free. Applying the inductive hypothesis to the graph $G'$, there exists an IO-code $S'$ of $G'$ such that $|S'| \le \left( \frac{2\Delta - 1}{2\Delta} \right)n$. We claim that $S'$ is also an IO-code of $G$. Suppose, to the contrary, that $S'$ is not an IO-code of $G$. In this case, at least one of $u$ and $v$ must be in the set $S'$. Renaming $u$ and $v$ if necessary, we may assume that $u \in S$ and that the vertex~$v$ is not identified by $S'$ in $G$. Let $x$ be the vertex in $G$ such that $v$ and $x$ are identified by $S'$ in $G'$ but are not identified by $S'$ in $G$. We note that the vertex $x$ is neighbor of $u$ different from~$v$. Since $S'$ is an IO-code of $G'$, in order for $S'$ to totally dominate the vertex $v$ in $G'$, there exists a vertex $w \in S'$ such that $w \ne u$ and $vw \in E(G)$. Since the set $S'$ does not identify the pair $v$ and $x$ in the graph $G'$, we infer that $w \ne x$ and $wx \in E(G)$. This implies that $G$ has a $4$-cycle, namely, $uvwxu$, a contradiction. Therefore, $S'$ is also an IO-code of $G$, and so $\ioc(G) \le |S'| \le \left( \frac{2\Delta - 1}{2\Delta} \right)n$.

Hence we may assume that $G'$ has open twins, for otherwise the desired result follows. More generally, we may assume that for every cycle edge $f$ of $G$, the graph $G - f$ has open twins. Since $G'$ has open twins, at least one of $u$ and $v$ is an open twin in $G'$. Renaming vertices if necessary, we may assume that $v$ is an open twin in $G'$. Let $x$ be an open twin of $v$ in $G'$. Since $G'$ has no $4$-cycles, we infer that the open twins $x$ and $v$ must be of degree~$1$ in $G'$. Let $w$ be the common neighbor of $x$ and $v$. Thus, $w$ is a support vertex in $G'$ with $v$ and $x$ as leaf neighbors. Since $e$ is a cycle edge of $G$, we note that the vertex $w$ belongs to the cycle $C$. Thus, $\deg_G(w) \ge 3$, $\deg_{G}(x) = 1$, $\deg_{G'}(v) = 1$ and $\deg_G(v) = 2$. Moreover, we cannot have $x = u$, or else, the cycle $C$ is $uvwu$. In other words, $C$ is a $3$-cycle in $G$. However, since $\deg_G(w) \ge 3$, it implies that removing the cycle edge $uw$ from $G$ does not produce open twins which contradicts our earlier assumption.

Let $f = vw \in E(G)$ and let $G'' = G - f$. We note that the edge $f$ belongs to the cycle $C$. By our earlier assumptions, the removal of the cycle edge $f$ creates open twins. By our earlier observations, open twins in $G''$ must be of degree~$1$ in $G''$. Since both $u$ and $w$ have degree at least~$2$ in $G''$, we therefore infer that the vertex $v$ is an open twin in $G''$. Let $y$ be an open twin of $v$ in $G''$. Analogously as in the previous case when $v$ is an open twin in the graph $G'$, we infer that the open twins $v$ and $y$ must be of degree~$1$ in $G''$. Thus, the vertex $u$ is a support vertex in $G''$ with $v$ and $y$ as leaf neighbors.

Continuing the same analysis for every edge of $C$, we infer that, if $C \colon v_0v_1v_2 \ldots v_p v_0$, then renaming vertices of $C$ if necessary, the vertex $v_i$ is a support vertex in $G$ for each odd $i \in [p]$ and the vertex $v_i$ has degree~$2$ in $G$ for each even $i \in [p]$. Thus, we note $p$ is an odd integer and, since $G$ contains no $4$-cycles, $p$ must be at least~$5$. Therefore, $C$ is a cycle of even length~$p + 1 \ge 6$. Let $x_i$ be a leaf neighbor of the vertex $v_i$ for each odd $i \in [p]$ (see Figure~\ref{fig:thm:main1}). We note that for every odd $i \in [p]$, the vertex $v_i$ has two neighbors on the cycle $C$ (and these neighbors are distinct from the leaf $x_i$), and so $v_i$ has degree at least~$3$ in $G$. Moreover, $m \ge p + p/2 = 3p/2 \ge 9$.

\vskip 0.25 cm 
\begin{figure}[htb]
\begin{center}
\begin{tikzpicture}[scale=.85,style=thick,x=0.85cm,y=0.85cm]
\def\vr{2.25pt}
\def\vrs{0.5pt}
\path (0,4.5) coordinate (d0);
\path (1,5) coordinate (d1);
\path (15.5,4.5) coordinate (d2);
\path (14.5,5) coordinate (d3);
\path (1,4) coordinate (u5);
\path (2.5,3.25) coordinate (v2);
\path (4.5,3.25) coordinate (v3);
\path (3.5,3) coordinate (v4);
\path (3.5,4) coordinate (v5);
\path (6,4) coordinate (w5);
\path (7.5,3.25) coordinate (x2);
\path (9.5,3.25) coordinate (x3);
\path (8.5,3) coordinate (x4);
\path (8.5,4) coordinate (x5);
\path (10,4) coordinate (p1);
\path (10.25,4) coordinate (p2);
\path (10.5,4) coordinate (p3);
\path (9.5,4) coordinate (q1);
\path (11.5,4) coordinate (q2);
\path (12,4) coordinate (y5);
\path (13.5,3.25) coordinate (z2);
\path (15.5,3.25) coordinate (z3);
\path (14.5,3) coordinate (z4);
\path (14.5,4) coordinate (z5);

\draw (v4)--(v5);
\draw (v3)--(v5);
\draw (v2)--(v5);
\draw (x4)--(x5);
\draw (x3)--(x5);
\draw (x2)--(x5);
\draw (z4)--(z5);
\draw (z3)--(z5);
\draw (z2)--(z5);
\draw (u5)--(v5)--(w5)--(x5);
\draw (y5)--(z5);
\draw (x5)--(q1);
\draw (y5)--(q2);
\draw (d1)--(d3);
\draw (d0) to[out=90,in=180, distance=0.5cm] (d1);
\draw (d0) to[out=270,in=180, distance=0.5cm] (u5);
\draw (d2) to[out=270,in=0, distance=0.5cm] (z5);
\draw (d2) to[out=90,in=0, distance=0.5cm] (d3);
\draw (u5) [fill=black] circle (\vr);
\draw (v4) [fill=black] circle (\vr);
\draw (v5) [fill=black] circle (\vr);
\draw (w5) [fill=black] circle (\vr);
\draw (x4) [fill=black] circle (\vr);
\draw (x5) [fill=black] circle (\vr);
\draw (p1) [fill=black] circle (\vrs);
\draw (p2) [fill=black] circle (\vrs);
\draw (p3) [fill=black] circle (\vrs);
\draw (y5) [fill=black] circle (\vr);
\draw (z4) [fill=black] circle (\vr);
\draw (z5) [fill=black] circle (\vr);
\draw[anchor = south] (u5) node {$v_0$};
\draw[anchor = south] (v5) node {$v_1$};
\draw[anchor = south] (w5) node {$v_2$};
\draw[anchor = south] (x5) node {$v_3$};
\draw[anchor = south] (y5) node {$v_{p-1}$};
\draw[anchor = south] (z5) node {$v_p$};
\draw (7.25,4.5) node {{\small $C$}};
\draw[anchor = north] (v4) node {$x_1$};
\draw[anchor = north] (x4) node {$x_3$};
\draw[anchor = north] (z4) node {$x_p$};
\end{tikzpicture}
\caption{A possible structure of the graph $G$ in the proof of Theorem~\ref{thm:main1}.}
\label{fig:thm:main1}
\end{center}
\end{figure}
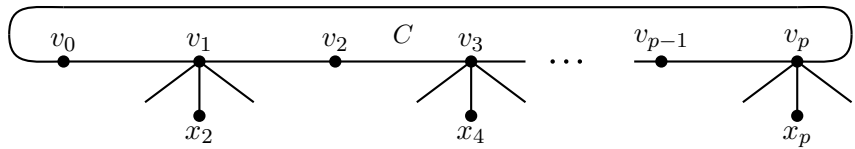

We now consider the graph $G_0 = G - v_0$ obtained by deleting the vertex $v_0$ from $G$. The graph $G_0$ is a connected graph satisfying $\Delta(G_0) \le \Delta(G) \le \Delta$. Since $G$ has no $4$-cycles, the graph $G_0$ also has no $4$-cycles. Recall that the neighbors $v_1$ and $v_p$ of $v_0$ are of degree at least~$3$ in $G$ and are support vertices with leaf neighbors $x_1$ and $x_p$, respectively. Since the graph $G$ is open twin-free, we therefore infer that the graph $G_0$ has no open twins. Let $G_0$ have size~$m_0$, and so $m_0 = m(G_0) = m - 2$. By our earlier observations, $m \ge 9$, and so $m_0 \ge 7$. Moreover since $p \ge 5$, the structure of the graph $G$ implies that $G_0 \ncong T_\Delta$ (noting that $x_1v_1v_2v_3x_3$ is a path in $G_0$ where $x_1$ and $x_3$ are leaves in $G_0$, the vertex $v_2$ has degree~$2$ in $G_0$, and $v_3$ is a vertex of degree at least~$3$ in $G_0$). Hence, by our induction hypothesis, there exists an IO-code $S_0$ of $G_0$ such that $|S_0| \le \left( \frac{2\Delta - 1}{2\Delta} \right)(n - 1) < \left( \frac{2\Delta - 1}{2\Delta} \right) n$.

We show next that $S_0$ is also an IO-code of $G$. Since $v_1$ and $v_p$ are support vertices in $G_0$, we note that $v_1, v_p \in S_0$. To show that $S_0$ is an IO-code of $G$, we only need to show that $S_0$ identifies the vertex $v_0$ from every other vertex of $G$. However, since $G$ has no $4$-cycles, this is indeed true, as $v_0$ is the only vertex of $G$ with $N_G(v_0) \cap S_0 = \{v_1, v_p\}$. Hence, the result is true in this case as well. This completes the proof of Theorem~\ref{thm:main1}.
\end{proof}

\section{Tight examples}
\label{S:construct}

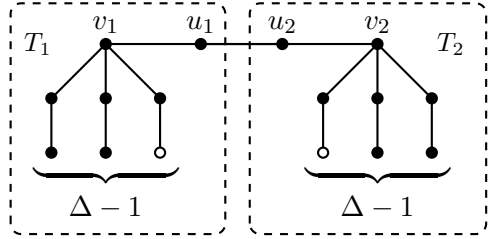
\begin{figure}[ht!]
\begin{center}
\begin{tikzpicture}[scale=.85,style=thick,x=0.85cm,y=0.85cm]
\def\vr{2.25pt}
\path (0,0) coordinate (u1);
\path (0,1) coordinate (v1);
\path (1,0) coordinate (u2);
\path (1,1) coordinate (v2);
\path (2,0) coordinate (u3);
\path (2,1) coordinate (v3);
\path (2.75,2) coordinate (u4);
\path (4.25,2) coordinate (v4);
\path (1,2) coordinate (v);
\path (5,0) coordinate (u5);
\path (5,1) coordinate (v5);
\path (6,0) coordinate (u6);
\path (6,1) coordinate (v6);
\path (7,0) coordinate (u7);
\path (7,1) coordinate (v7);
\path (6,2) coordinate (u);
\draw (v)--(v1)--(u1);
\draw (v)--(v2)--(u2);
\draw (v)--(v3)--(u3);
\draw (v)--(v4)--(u4);
\draw (u)--(v5)--(u5);
\draw (u)--(v6)--(u6);
\draw (u)--(v7)--(u7);
\draw (u)--(v4);
\draw (1,-0.25) node {$\underbrace{\phantom{1111111111}}$};
\draw (1,-1) node {$\Delta - 1$};
\draw (6,-0.25) node {$\underbrace{\phantom{1111111111}}$};
\draw (6,-1) node {$\Delta - 1$};
\draw [style=dashed,rounded corners] (-0.75,-1.5) rectangle (3.2,2.75);
\draw [style=dashed,rounded corners] (3.65,-1.5) rectangle (8,2.75);
\draw (-0.25,2) node {{\small $T_1$}};
\draw (7.35,2) node {{\small $T_2$}};
\draw (v1) [fill=black] circle (\vr);
\draw (v2) [fill=black] circle (\vr);
\draw (v3) [fill=black] circle (\vr);
\draw (v4) [fill=black] circle (\vr);
\draw (v5) [fill=black] circle (\vr);
\draw (v6) [fill=black] circle (\vr);
\draw (v7) [fill=black] circle (\vr);
\draw (u1) [fill=black] circle (\vr);
\draw (u2) [fill=black] circle (\vr);
\draw (u3) [fill=white] circle (\vr);
\draw (u4) [fill=black] circle (\vr);
\draw (u5) [fill=white] circle (\vr);
\draw (u6) [fill=black] circle (\vr);
\draw (u7) [fill=black] circle (\vr);
\draw (u) [fill=black] circle (\vr);
\draw (v) [fill=black] circle (\vr);
\draw[anchor = south] (v) node {$v_1$};
\draw[anchor = south] (u) node {$v_2$};
\draw[anchor = south] (u4) node {$u_1$};
\draw[anchor = south] (v4) node {$u_2$};
\end{tikzpicture}
\caption{A tree $T$ of order~$n$ satisfying $\ioc(T) = \left( \frac{2\Delta - 1}{2\Delta} \right) n$.}
\label{fig:bounds1}
\end{center}
\end{figure}

For $\Delta \ge 3$ a fixed integer, let $T_1$ and $T_2$ be two vertex-disjoint copies of the reduced subdivided star $T_{\Delta}^*$. Let $v_i$ be the central vertex in $T_i$ (of degree~$\Delta$) and let $u_i$ be the leaf neighbor of $v_i$ for $i \in [2]$. Let $T_{1,2}$ be the tree obtained from the union of $T_1$ and $T_2$ by adding the edge $u_1u_2$. The resulting tree $T$ is illustrated in Figure~\ref{fig:bounds1} and satisfies $\ioc(T) = \left( \frac{2\Delta - 1}{2\Delta} \right) n$. These examples show that the upper bound in Theorem~\ref{thm:main1} is best possible for every fixed value of $\Delta \ge 3$.

\begin{figure}[ht!]
\begin{center}
\begin{tikzpicture}[scale=.85,style=thick,x=0.85cm,y=0.85cm]
\def\vr{2.25pt}
\def\vrs{0.5pt}
\path (0,4.5) coordinate (d0);
\path (1,5) coordinate (d1);
\path (15.5,4.5) coordinate (d2);
\path (14.5,5) coordinate (d3);
\path (0,2) coordinate (u0);
\path (1,0) coordinate (u1);
\path (1,1) coordinate (u2);
\path (1,2) coordinate (u3);
\path (1,3) coordinate (u4);
\path (1,4) coordinate (u5);
\path (2.5,2) coordinate (v0);
\path (3.5,0) coordinate (v1);
\path (3.5,1) coordinate (v2);
\path (3.5,2) coordinate (v3);
\path (3.5,3) coordinate (v4);
\path (3.5,4) coordinate (v5);
\path (5,2) coordinate (w0);
\path (6,0) coordinate (w1);
\path (6,1) coordinate (w2);
\path (6,2) coordinate (w3);
\path (6,3) coordinate (w4);
\path (6,4) coordinate (w5);
\path (7.5,2) coordinate (x0);
\path (8.5,0) coordinate (x1);
\path (8.5,1) coordinate (x2);
\path (8.5,2) coordinate (x3);
\path (8.5,3) coordinate (x4);
\path (8.5,4) coordinate (x5);
\path (9.5,4) coordinate (p1);
\path (9.75,4) coordinate (p2);
\path (10,4) coordinate (p3);
\path (9,4) coordinate (q1);
\path (10.5,4) coordinate (q2);
\path (10.5,2) coordinate (y0);
\path (11.5,0) coordinate (y1);
\path (11.5,1) coordinate (y2);
\path (11.5,2) coordinate (y3);
\path (11.5,3) coordinate (y4);
\path (11.5,4) coordinate (y5);
\path (13.5,2) coordinate (z0);
\path (14.5,0) coordinate (z1);
\path (14.5,1) coordinate (z2);
\path (14.5,2) coordinate (z3);
\path (14.5,3) coordinate (z4);
\path (14.5,4) coordinate (z5);

\draw (u1)--(u2)--(u3)--(u4)--(u5);
\draw (u3)--(u0);
\draw (v1)--(v2)--(v3)--(v4)--(v5);
\draw (v3)--(v0);
\draw (w1)--(w2)--(w3)--(w4)--(w5);
\draw (w3)--(w0);
\draw (x1)--(x2)--(x3)--(x4)--(x5);
\draw (x3)--(x0);
\draw (y1)--(y2)--(y3)--(y4)--(y5);
\draw (y3)--(y0);
\draw (z1)--(z2)--(z3)--(z4)--(z5);
\draw (z3)--(z0);
\draw (u5)--(v5)--(w5)--(x5);
\draw (y5)--(z5);
\draw (x5)--(q1);
\draw (y5)--(q2);
\draw (d1)--(d3);
\draw (d0) to[out=90,in=180, distance=0.5cm] (d1);
\draw (d0) to[out=270,in=180, distance=0.5cm] (u5);
\draw (d2) to[out=270,in=0, distance=0.5cm] (z5);
\draw (d2) to[out=90,in=0, distance=0.5cm] (d3);
\draw (u0) [fill=white] circle (\vr);
\draw (u1) [fill=black] circle (\vr);
\draw (u2) [fill=black] circle (\vr);
\draw (u3) [fill=black] circle (\vr);
\draw (u4) [fill=black] circle (\vr);
\draw (u5) [fill=black] circle (\vr);
\draw (v0) [fill=white] circle (\vr);
\draw (v1) [fill=black] circle (\vr);
\draw (v2) [fill=black] circle (\vr);
\draw (v3) [fill=black] circle (\vr);
\draw (v4) [fill=black] circle (\vr);
\draw (v5) [fill=black] circle (\vr);
\draw (w0) [fill=white] circle (\vr);
\draw (w1) [fill=black] circle (\vr);
\draw (w2) [fill=black] circle (\vr);
\draw (w3) [fill=black] circle (\vr);
\draw (w4) [fill=black] circle (\vr);
\draw (w5) [fill=black] circle (\vr);
\draw (x0) [fill=white] circle (\vr);
\draw (x1) [fill=black] circle (\vr);
\draw (x2) [fill=black] circle (\vr);
\draw (x3) [fill=black] circle (\vr);
\draw (x4) [fill=black] circle (\vr);
\draw (x5) [fill=black] circle (\vr);
\draw (p1) [fill=black] circle (\vrs);
\draw (p2) [fill=black] circle (\vrs);
\draw (p3) [fill=black] circle (\vrs);
\draw (y0) [fill=white] circle (\vr);
\draw (y1) [fill=black] circle (\vr);
\draw (y2) [fill=black] circle (\vr);
\draw (y3) [fill=black] circle (\vr);
\draw (y4) [fill=black] circle (\vr);
\draw (y5) [fill=black] circle (\vr);
\draw (z0) [fill=white] circle (\vr);
\draw (z1) [fill=black] circle (\vr);
\draw (z2) [fill=black] circle (\vr);
\draw (z3) [fill=black] circle (\vr);
\draw (z4) [fill=black] circle (\vr);
\draw (z5) [fill=black] circle (\vr);
\draw[anchor = south] (u5) node {$u_1$};
\draw[anchor = south] (v5) node {$u_2$};
\draw[anchor = south] (w5) node {$u_3$};
\draw[anchor = south] (x5) node {$u_4$};
\draw[anchor = south] (y5) node {$u_{p-1}$};
\draw[anchor = south] (z5) node {$u_p$};
\draw (7.25,4.5) node {{\small $C$}};
\draw[anchor = south] (u0) node {$z_1$};
\draw[anchor = west] (u4) node {$v_1$};
\draw[anchor = west] (u3) node {$w_1$};
\draw[anchor = west] (u2) node {$x_1$};
\draw[anchor = west] (u1) node {$y_1$};
\draw[anchor = south] (v0) node {$z_2$};
\draw[anchor = west] (v4) node {$v_2$};
\draw[anchor = west] (v3) node {$w_2$};
\draw[anchor = west] (v2) node {$x_2$};
\draw[anchor = west] (v1) node {$y_2$};
\draw[anchor = south] (w0) node {$z_3$};
\draw[anchor = west] (w4) node {$v_3$};
\draw[anchor = west] (w3) node {$w_3$};
\draw[anchor = west] (w2) node {$x_3$};
\draw[anchor = west] (w1) node {$y_3$};
\draw[anchor = south] (x0) node {$z_4$};
\draw[anchor = west] (x4) node {$v_4$};
\draw[anchor = west] (x3) node {$w_4$};
\draw[anchor = west] (x2) node {$x_4$};
\draw[anchor = west] (x1) node {$y_4$};
\draw[anchor = south] (y0) node {$z_{p-1}$};
\draw[anchor = west] (y4) node {$v_{p-1}$};
\draw[anchor = west] (y3) node {$w_{p-1}$};
\draw[anchor = west] (y2) node {$x_{p-1}$};
\draw[anchor = west] (y1) node {$y_{p-1}$};
\draw[anchor = south] (z0) node {$z_p$};
\draw[anchor = west] (z4) node {$v_p$};
\draw[anchor = west] (z3) node {$w_p$};
\draw[anchor = west] (z2) node {$x_p$};
\draw[anchor = west] (z1) node {$y_p$};
\end{tikzpicture}
\caption{A subcubic graph $G_p$ of order~$n$ satisfying $\ioc(G_p) = \frac{5}{6}n$.}
\label{fig:subcubic}
\end{center}
\end{figure}
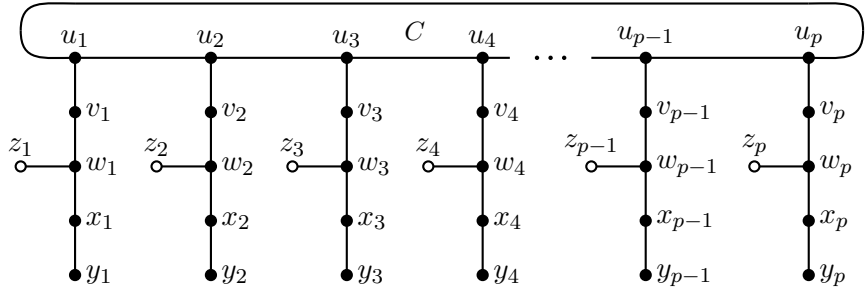

A \emph{subcubic graph} is a graph with maximum degree at most~$3$. In the special case when $\Delta = 3$, by Theorem~\ref{thm:main1} if $G$ is a subcubic graph of order~$n \ge 5$ that is open twin-free and contains no $4$-cycles, then $\ioc(G) \le \frac{5}{6}n$. We construct next a family of subcubic graphs $G$ of arbitrarily large orders~$n$ that are open twin-free and contain no $4$-cycles satisfying $\ioc(G) = \frac{5}{6}n$. Let $C \colon u_1u_2 \ldots u_pu_1$ be a cycle of length~$p$ where $p \ge 3$ and $p \ne 4$. For each vertex $u_i$ on the cycle, we add a vertex disjoint copy, $T_i$ say, of a reduced subdivided star $T_{3}^*$ and identify one of the leaves at distance~$2$ from the central vertex of the reduced subdivided star with the vertex~$u_i$ for all $i \in [p]$. Let $V(T_i) = \{u_i,v_i,w_i,x_i,y_i,z_i\}$ where $u_iv_iw_ix_iy_i$ is a path and where $z_i$ is the leaf neighbor of $w_i$ in $T_i$. The resulting subcubic graph $G_p$ is illustrated in Figure~\ref{fig:subcubic}.

\begin{proposition}
\label{prop:subcubic}
Let $p \ge 3$ and $p \ne 4$. Moreover, let $G_p$ be a subcubic graph of order $n$. Then $G_p$ is open twin-free, contains no $4$-cycles and satisfies $\ioc(G_p) = \frac{5}{6}n$.
\end{proposition}

\begin{proof}
 Let $S$ be an arbitrary IO-code in $G_p$. We show that $|S \cap V(T_i)| \ge 5$ for all $i \in [p]$. Since $S$ is a TD-set, the set $S$ contains the support vertices $w_i$ and $x_i$. In order to identify $x_i$ and $z_i$, the vertex $y_i$ belongs to the set $S$.  In order to identify $v_i$ and $z_i$, the vertex $u_i$ belongs to the set $S$. In order to identify $w_i$ and $y_i$, the set $S$ contains at least one of $v_i$ and $z_i$. Hence, $|S \cap V(T_i)| \ge 5$ for all $i \in [p]$. Since $S$ is an arbitrary IO-code in $G_p$ and $|S| \ge \frac{5}{6}n$, this implies that $\ioc(T) \ge \frac{5}{6}n$. Since the set $S^* = V(G_p) \setminus \{z_1,z_2,\ldots,z_p\}$ is an IO-code in $G_p$, we have $\ioc(T) \le \frac{5}{6}n$. Consequently, $\ioc(T) = \frac{5}{6}n$.
 \end{proof}

\medskip
We remark that deleting the cycle edge $u_1u_p$ in the construction of the subcubic graph $G_p$ yields a tree $T$ of order~$n$ that is open twin-free and satisfies $\ioc(T) = \frac{5}{6}n$. By Proposition~\ref{prop:subcubic}, the upper bound in Theorem~\ref{thm:main1} is tight for $\Delta = 3$ in the strong sense that there exist connected subcubic graphs of arbitrary large order~$n$ that contain no $4$-cycles, are open twin-free, and that achieve the upper bound $\ioc(G) = \left( \frac{2\Delta - 1}{2\Delta} \right) n$.

\medskip

For a fixed $\Delta \ge 4$ however, we believe that there are no tight examples of the upper bound in Theorem~\ref{thm:main1}. Nevertheless, as we shall show below, the said upper bound in Theorem~\ref{thm:main1} is still nearly tight. 
To demonstrate this fact, we construct a family of graphs, denoted $G_{p,\Delta}$ where $p \ge 3$ with $p \ne 4$ and $\Delta \ge 4$, of arbitrarily large orders~$n$ and with maximum degree $\Delta$ that are open twin-free, contain no $4$-cycles and whose IO-code numbers attain the value $\left( \frac{2\Delta-4}{2\Delta-3} \right) n$. As we show below, this implies that for this family of graphs $G_{p, \Delta}$, the parameter $\ioc(G_{p,\Delta})$ falls short of the upper bound in Theorem~\ref{thm:main1} by $\frac{3p}{2\Delta}$, which for fixed $p$ asymptotically tends to zero with large $\Delta$. 

Let $\Delta \ge 4$ be a fixed integer and let $C \colon u_1u_2 \ldots u_pu_1$ be a cycle of length~$p$, where $p \ge 3$ and $p \ne 4$. For each vertex $u_i$ on the cycle, we introduce a vertex-disjoint copy, $T_i$ say, of a subdivided star $T_{\Delta-2}$. Let $V(T_i) = \{w_i,v_i^1, x_i^1, v_i^2, x_i^2, v_i^3, x_i^3, \ldots , v_i^{\Delta-2}, x_i^{\Delta-2}\}$, where $w_i$ is the central vertex of $T_i$, each $v_i^j$ for $1 \le j \le \Delta-2$ is a support vertex of $T_i$ at distance~$1$ from $w_i$, and each $x_i^j$ is a leaf of $T_i$, adjacent to $v_i^j$, and at distance~$2$ from $w_i$. Let the graph constructed be called $G_{p,\Delta}$. An example of such a graph is illustrated in Figure~\ref{fig:Delta>=4}. Notice that $n_{p,\Delta} = n(G_{p,\Delta}) = p (2\Delta-3)$. Also notice that the vertices $u_i$ for $i \in [p]$ are of degree~$\Delta$ each.

\begin{figure}[t!]
\begin{center}
\begin{tikzpicture}[scale=0.85,style=thick,x=0.85cm,y=0.85cm]
\def\vr{2.25pt}
\def\vrs{0.5pt}


\path (-2,2.5) coordinate (d0);
\path (-1,3) coordinate (d1);
\path (16,2.5) coordinate (d2);
\path (15,3) coordinate (d3);

\path (-1,2) coordinate (u5);
\path (3,2) coordinate (v5);
\path (7,2) coordinate (w5);
\path (11,2) coordinate (y5);
\path (15,2) coordinate (z5);

\path (8.75,2) coordinate (p1);
\path (9,2) coordinate (p2);
\path (9.25,2) coordinate (p3);

\draw (u5)--(v5)--(w5)--(8.3,2);
\draw (9.7,2)--(y5)--(z5);

\draw (u5) [fill=white] circle (\vr);
\draw (v5) [fill=white] circle (\vr);
\draw (w5) [fill=white] circle (\vr);
\draw (p1) [fill=black] circle (\vrs);
\draw (p2) [fill=black] circle (\vrs);
\draw (p3) [fill=black] circle (\vrs);
\draw (y5) [fill=white] circle (\vr);
\draw (z5) [fill=white] circle (\vr);


\draw[anchor = south] (u5) node {$u_1$};
\draw[anchor = south] (v5) node {$u_2$};
\draw[anchor = south] (w5) node {$u_3$};
\draw[anchor = south] (y5) node {$u_{p-1}$};
\draw[anchor = south] (z5) node {$u_p$};
\draw (9.25,2.5) node {{\small $C$}};
\draw (d1)--(d3);
\draw (d0) to[out=90,in=180, distance=0.5cm] (d1);
\draw (d0) to[out=270,in=180, distance=0.5cm] (u5);
\draw (d2) to[out=270,in=0, distance=0.5cm] (z5);
\draw (d2) to[out=90,in=0, distance=0.5cm] (d3);

%
\path (-2.5,0) coordinate (w1);
\path (-2.5,1) coordinate (x1);
\path (-1.5,0) coordinate (w2);
\path (-1.5,1) coordinate (x2);
\path (-0.5,0) coordinate (w3);
\path (-0.5,1) coordinate (x3);
\path (0.5,0) coordinate (w4);
\path (0.5,1) coordinate (x4);
\path (-1,2) coordinate (w);
\path (-0.15,0) coordinate (q1);
\path (0,0) coordinate (q2);
\path (0.15,0) coordinate (q3);

\draw (w)--(x1)--(w1);
\draw (w)--(x2)--(w2);
\draw (w)--(x3)--(w3);
\draw (w)--(x4)--(w4);
\draw (-1,-0.7) node {$\underbrace{\phantom{1111111111111}}$};
\draw (-1,-1.45) node {$\Delta-2$};
\draw (w1) [fill=white] circle (\vr);
\draw (w2) [fill=black] circle (\vr);
\draw (w3) [fill=black] circle (\vr);
\draw (w4) [fill=black] circle (\vr);
\draw (x1) [fill=black] circle (\vr);
\draw (x2) [fill=black] circle (\vr);
\draw (x3) [fill=black] circle (\vr);
\draw (x4) [fill=black] circle (\vr);
\draw (w) [fill=black] circle (\vr);
\draw (p1) [fill=black] circle (\vrs);
\draw (p2) [fill=black] circle (\vrs);
\draw (p3) [fill=black] circle (\vrs);
\draw (q1) [fill=black] circle (\vrs);
\draw (q2) [fill=black] circle (\vrs);
\draw (q3) [fill=black] circle (\vrs);
{\footnotesize
\draw[anchor = east] (x1) node {$v^1_1$};
\draw[anchor = east] (x2) node {$v^2_1$};
\draw[anchor = east] (x3) node {$v^3_1$};
\node () at (0.6,1.5) {$v^{\Delta-2}_1$};

\draw[anchor = north] (w1) node {$x^1_1$};
\draw[anchor = north] (w2) node {$x^2_1$};
\draw[anchor = north] (w3) node {$x^3_1$};
\draw[anchor = north] (w4) node {$x^{\Delta-2}_1$};}


\path (1.5,0) coordinate (w1');
\path (1.5,1) coordinate (x1');
\path (2.5,0) coordinate (w2');
\path (2.5,1) coordinate (x2');
\path (3.5,0) coordinate (w3');
\path (3.5,1) coordinate (x3');
\path (4.5,0) coordinate (w4');
\path (4.5,1) coordinate (x4');
\path (3,2) coordinate (w');
\path (3.85,0) coordinate (q1');
\path (4,0) coordinate (q2');
\path (4.15,0) coordinate (q3');

\draw (w')--(x1')--(w1');
\draw (w')--(x2')--(w2');
\draw (w')--(x3')--(w3');
\draw (w')--(x4')--(w4');
\draw (3,-0.7) node {$\underbrace{\phantom{1111111111111}}$};
\draw (3,-1.45) node {$\Delta-2$};
\draw (w1') [fill=white] circle (\vr);
\draw (w2') [fill=black] circle (\vr);
\draw (w3') [fill=black] circle (\vr);
\draw (w4') [fill=black] circle (\vr);
\draw (x1') [fill=black] circle (\vr);
\draw (x2') [fill=black] circle (\vr);
\draw (x3') [fill=black] circle (\vr);
\draw (x4') [fill=black] circle (\vr);
\draw (w') [fill=black] circle (\vr);
\draw (q1') [fill=black] circle (\vrs);
\draw (q2') [fill=black] circle (\vrs);
\draw (q3') [fill=black] circle (\vrs);

{\footnotesize
\draw[anchor = east] (x1') node {$v^1_2$};
\draw[anchor = east] (x2') node {$v^2_2$};
\draw[anchor = east] (x3') node {$v^3_2$};
\node () at (4.6,1.5) {$v^{\Delta-2}_2$};

\draw[anchor = north] (w1') node {$x^1_2$};
\draw[anchor = north] (w2') node {$x^2_2$};
\draw[anchor = north] (w3') node {$x^3_2$};
\draw[anchor = north] (w4') node {$x^{\Delta-2}_2$};}


\path (5.5,0) coordinate (w1'');
\path (5.5,1) coordinate (x1'');
\path (6.5,0) coordinate (w2'');
\path (6.5,1) coordinate (x2'');
\path (7.5,0) coordinate (w3'');
\path (7.5,1) coordinate (x3'');
\path (8.5,0) coordinate (w4'');
\path (8.5,1) coordinate (x4'');
\path (7,2) coordinate (w'');
\path (7.85,0) coordinate (q1'');
\path (8,0) coordinate (q2'');
\path (8.15,0) coordinate (q3'');

\draw (w'')--(x1'')--(w1'');
\draw (w'')--(x2'')--(w2'');
\draw (w'')--(x3'')--(w3'');
\draw (w'')--(x4'')--(w4'');
\draw (7,-0.7) node {$\underbrace{\phantom{1111111111111}}$};
\draw (7,-1.45) node {$\Delta-2$};
\draw (w1'') [fill=white] circle (\vr);
\draw (w2'') [fill=black] circle (\vr);
\draw (w3'') [fill=black] circle (\vr);
\draw (w4'') [fill=black] circle (\vr);
\draw (x1'') [fill=black] circle (\vr);
\draw (x2'') [fill=black] circle (\vr);
\draw (x3'') [fill=black] circle (\vr);
\draw (x4'') [fill=black] circle (\vr);
\draw (w'') [fill=black] circle (\vr);
\draw (q1'') [fill=black] circle (\vrs);
\draw (q2'') [fill=black] circle (\vrs);
\draw (q3'') [fill=black] circle (\vrs);
{\footnotesize
\draw[anchor = east] (x1'') node {$v^1_3$};
\draw[anchor = east] (x2'') node {$v^2_3$};
\draw[anchor = east] (x3'') node {$v^3_3$};
\node () at (8.6,1.5) {$v^{\Delta-2}_3$};

\draw[anchor = north] (w1'') node {$x^1_3$};
\draw[anchor = north] (w2'') node {$x^2_3$};
\draw[anchor = north] (w3'') node {$x^3_3$};
\draw[anchor = north] (w4'') node {$x^{\Delta-2}_3$};}


\path (9.5,0) coordinate (w1''');
\path (9.5,1) coordinate (x1''');
\path (10.5,0) coordinate (w2''');
\path (10.5,1) coordinate (x2''');
\path (11.5,0) coordinate (w3''');
\path (11.5,1) coordinate (x3''');
\path (12.5,0) coordinate (w4''');
\path (12.5,1) coordinate (x4''');
\path (11,2) coordinate (w''');
\path (11.85,0) coordinate (q1''');
\path (12,0) coordinate (q2''');
\path (12.15,0) coordinate (q3''');

\draw (w''')--(x1''')--(w1''');
\draw (w''')--(x2''')--(w2''');
\draw (w''')--(x3''')--(w3''');
\draw (w''')--(x4''')--(w4''');
\draw (11,-0.7) node {$\underbrace{\phantom{1111111111111}}$};
\draw (11,-1.45) node {$\Delta-2$};
\draw (w1''') [fill=white] circle (\vr);
\draw (w2''') [fill=black] circle (\vr);
\draw (w3''') [fill=black] circle (\vr);
\draw (w4''') [fill=black] circle (\vr);
\draw (x1''') [fill=black] circle (\vr);
\draw (x2''') [fill=black] circle (\vr);
\draw (x3''') [fill=black] circle (\vr);
\draw (x4''') [fill=black] circle (\vr);
\draw (w''') [fill=black] circle (\vr);
\draw (q1''') [fill=black] circle (\vrs);
\draw (q2''') [fill=black] circle (\vrs);
\draw (q3''') [fill=black] circle (\vrs);
{\footnotesize
\draw[anchor = east] (x1''') node {$v^1_4$};
\draw[anchor = east] (x2''') node {$v^2_4$};
\draw[anchor = east] (x3''') node {$v^3_4$};
\node () at (12.6,1.5) {$v^{\Delta-2}_4$};

\draw[anchor = north] (w1''') node {$x^1_4$};
\draw[anchor = north] (w2''') node {$x^2_4$};
\draw[anchor = north] (w3''') node {$x^3_4$};
\draw[anchor = north] (w4''') node {$x^{\Delta-2}_4$};}


\path (13.5,0) coordinate (w1'''');
\path (13.5,1) coordinate (x1'''');
\path (14.5,0) coordinate (w2'''');
\path (14.5,1) coordinate (x2'''');
\path (15.5,0) coordinate (w3'''');
\path (15.5,1) coordinate (x3'''');
\path (16.5,0) coordinate (w4'''');
\path (16.5,1) coordinate (x4'''');
\path (15,2) coordinate (w'''');
\path (15.85,0) coordinate (q1'''');
\path (16,0) coordinate (q2'''');
\path (16.15,0) coordinate (q3'''');

\draw (w'''')--(x1'''')--(w1'''');
\draw (w'''')--(x2'''')--(w2'''');
\draw (w'''')--(x3'''')--(w3'''');
\draw (w'''')--(x4'''')--(w4'''');
\draw (15,-0.7) node {$\underbrace{\phantom{1111111111111}}$};
\draw (15,-1.45) node {$\Delta-2$};
\draw (w1'''') [fill=white] circle (\vr);
\draw (w2'''') [fill=black] circle (\vr);
\draw (w3'''') [fill=black] circle (\vr);
\draw (w4'''') [fill=black] circle (\vr);
\draw (x1'''') [fill=black] circle (\vr);
\draw (x2'''') [fill=black] circle (\vr);
\draw (x3'''') [fill=black] circle (\vr);
\draw (x4'''') [fill=black] circle (\vr);
\draw (w'''') [fill=black] circle (\vr);
\draw (q1'''') [fill=black] circle (\vrs);
\draw (q2'''') [fill=black] circle (\vrs);
\draw (q3'''') [fill=black] circle (\vrs);
{\footnotesize
\draw[anchor = east] (x1'''') node {$v^1_p$};
\draw[anchor = east] (x2'''') node {$v^2_p$};
\draw[anchor = east] (x3'''') node {$v^3_p$};
\node () at (16.6,1.5) {$v^{\Delta-2}_p$};

\draw[anchor = north] (w1'''') node {$x^1_p$};
\draw[anchor = north] (w2'''') node {$x^2_p$};
\draw[anchor = north] (w3'''') node {$x^3_p$};
\draw[anchor = north] (w4'''') node {$x^{\Delta-2}_p$};}
\end{tikzpicture}
\caption{An open twin-free graph $G_{p,\Delta}$ of order~$n_{p,\Delta}$, with no $4$-cycles and with maximum degree $\Delta(G) \ge 4$ satisfying $\ioc(G_{p,\Delta}) = \left( \frac{2\Delta - 4}{2\Delta-3} \right) n_{p,\Delta}$}.
\label{fig:Delta>=4}
\end{center}
\end{figure}

\begin{proposition}
\label{prop:large_delta}
For $p \ge 3$ with $p \ne 4$ and $\Delta \ge 4$, if the graph $G_{p,\Delta}$ has order $n_{p,\Delta}$, then $G_{p,\Delta}$ is open twin-free, contains no $4$-cycles and satisfies $\ioc(G_{p,\Delta}) = \left( \frac{2\Delta-4}{2\Delta-3} \right) n_{p, \Delta}$.
\end{proposition}

\begin{proof}
Notice that the edges $u_iu_{i+1}$, for $i \in [p-1]$, of the cycle $C$ do not transform any of the leaves or support vertices of each $T_i$ into non-leaves or non-support vertices, respectively. Therefore by the exact same proof as for Proposition~\ref{prop:prop1}(a), any IO-code of $G_{p,\Delta}$ must contain at least all vertices but one from each $T_i$. Thus, we must have $\ioc(G_{p,\Delta}) \ge \sum_{1 \le i \le p} \left( n(T_i)-1 \right) = \sum_{1 \le i \le p} \left( \frac{2(\Delta-2)}{2(\Delta-2)+1} \right) n(T_i) = \left( \frac{2\Delta-4}{2\Delta-3} \right) n_{p,\Delta}$. Moreover, it can be checked that the set $S = \cup_{1 \le i \le p} (V(T_i) \setminus \{x_i^1\})$ is an IO-code of $G_{p, \Delta}$. Thus, we also have $\ioc(G_{p,\Delta}) \le |S| = p(2\Delta-4) = \left( \frac{2\Delta-4}{2\Delta-3} \right) n_{p,\Delta}$. This proves the result.
\end{proof}

Clearly, the parameter $\ioc(G_{p,\Delta})$ shown in Proposition~\ref{prop:large_delta} is smaller than the upper bound in Theorem~\ref{thm:main1} for the graphs $G_{p,\Delta}$. However, to compare this difference, we see that
\[
\left( \frac{2\Delta-1}{2\Delta} \right) n_{p,\Delta} - \left( \frac{2\Delta-4}{2\Delta-3} \right) n_{p,\Delta} = \frac{3}{2\Delta(2\Delta-3)} n_{p,\Delta} = \frac{3}{2\Delta(2\Delta-3)}p(2\Delta-3) = \frac{3p}{2\Delta}.
\]

In other words, irrespective of the length $p$ of the cycle $C$ (and hence, the order of the graph $G_{p,\Delta}$), we can approach the upper bound in Theorem~\ref{thm:main1} as close as necessary by taking large enough $\Delta \ge 4$.

\begin{rmk} The above construction can be slightly improved by adding more levels, in the following way. We start with a cycle $C_p$, where $p$ is a multiple of~$3$. To every other vertex of the cycle, we add $\Delta-2$ new neighbors, for each new neighbor, we further add $\Delta-1$ neighbors, and to each such leaf we attach a copy of $T_{\Delta-1}$. One can check that the IO-code number of the resulting graph is $\left(\frac{2\Delta-4}{2\Delta-3}+\frac{2}{f(\Delta)}\right)n$, where $f(\Delta) =  3(\Delta-1)^2(2\Delta-3)$.
\end{rmk}

As $f$ is a degree~$3$ polynomial, this is only a small improvement but significantly complicates the construction, and hence we do not give details. We did not manage to create examples with a higher ratio, such as for example $\left(\frac{2\Delta-3}{2\Delta-2}\right)n$, and leave this as an open problem.

\section{Concluding remarks}
\label{S:conclude}

Our main result, namely Theorem~\ref{thm:main1}, shows that for $\Delta \ge 3$ a fixed integer, if $G$ is a connected graph of order~$n \ge 5$ that is open twin-free and satisfies $\Delta(G) \le \Delta$, then $\ioc(G) \le \left( \frac{2\Delta - 1}{2\Delta} \right) n$, except in one exceptional case when $G$ is a subdivided star of maximum degree~$\Delta$, that is, if $G \cong  T_{\Delta}$. As shown in Proposition~\ref{prop:prop1}, if $G \cong T_{\Delta}^*$ is the reduced subdivided star of order~$n$, then $\ioc(G) = \left( \frac{2\Delta - 1}{2\Delta} \right) n$.

As remarked earlier, the upper bound in Theorem~\ref{thm:main1} is best possible when $\Delta = 3$, for arbitrarily large connected graphs. For $\Delta \ge 4$ however, it remains an open problem to determine if the upper bound in Theorem~\ref{thm:main1} is tight for arbitrarily large connected graphs. Nevertheless, for such values of $\Delta \ge 4$, we have presented examples of arbitrarily large open twin-free graphs of maximum degree at most $\Delta$ and without $4$-cycles whose IO-code numbers almost reach the upper bound in Theorem~\ref{thm:main1}. Yet, for $\Delta \ge 4$ there is room for small improvements, for example by constructing graphs whose IO-code numbers are even closer to the bound, such as $\left(\frac{2\Delta-3}{2\Delta-2}\right)n$ for example.

When we consider general graphs and thus allow $4$-cycles, the bound of Theorem~\ref{thm:main1} does not hold. In fact, as shown in~\cite{FoGhSh-21}, the infinite family of \emph{half-graphs} provides infinitely many connected bipartite graphs $G$ of order $n$ with $\ioc(G)=n$. By using these graphs as building blocks, one can build, for every fixed value of $\Delta$, arbitrarily large connected graphs $G$ of order~$n$ with $\ioc(G)=\left( \frac{2\Delta}{2\Delta+1} \right) n$ (similar to the graphs constructed in Proposition~\ref{prop:subcubic}). Start with a cycle $C_p$ of order $p$; for each vertex $v$ of $C_p$, consider a copy of the half-graph $H_\Delta$ and add an edge between $v$ and a vertex of degree~1 of $H_\Delta$. (Actually, $v$ can be made adjacent to any vertex of the half-graph that has degree strictly less than $\Delta$.) An example for $\Delta=3$ is provided in Figure~\ref{fig:cycle-halfgraphs}. Thus, one may ask what is the best possible bound of this form in the general case (when one excludes half-graphs): perhaps it is indeed $\left( \frac{2\Delta}{2\Delta+1} \right) n$?

\begin{figure}[ht!]
\begin{center}
\begin{tikzpicture}[scale=.85,style=thick,x=0.85cm,y=0.85cm]
\def\vr{2.25pt}
\def\vrs{0.5pt}
\path (7.5,4.5) coordinate (d0);
\path (8.5,5) coordinate (d1);
\path (15.5,4.5) coordinate (d2);
\path (14.5,5) coordinate (d3);
%

\path (8.5,4) coordinate (x0);
\path (8.5,3) coordinate (x1);
\path (8.5,2) coordinate (x2);
\path (8.5,1) coordinate (x3);
\path (7.5,3) coordinate (x4);
\path (7.5,2) coordinate (x5);
\path (7.5,1) coordinate (x6);
\path (9.5,4) coordinate (p1);
\path (9.75,4) coordinate (p2);
\path (10,4) coordinate (p3);
\path (9,4) coordinate (q1);
\path (10.5,4) coordinate (q2);
\path (11.5,4) coordinate (y0);
\path (11.5,3) coordinate (y1);
\path (11.5,2) coordinate (y2);
\path (11.5,1) coordinate (y3);
\path (10.5,3) coordinate (y4);
\path (10.5,2) coordinate (y5);
\path (10.5,1) coordinate (y6);
\path (14.5,4) coordinate (z0);
\path (14.5,3) coordinate (z1);
\path (14.5,2) coordinate (z2);
\path (14.5,1) coordinate (z3);
\path (13.5,3) coordinate (z4);
\path (13.5,2) coordinate (z5);
\path (13.5,1) coordinate (z6);

%
\draw (x0)--(x1)--(x4)--(x2)--(x5)--(x3);
\draw (x4)--(x3)--(x6);
\draw (y0)--(y1)--(y4)--(y2)--(y5)--(y3);
\draw (y4)--(y3)--(y6);
\draw (z0)--(z1)--(z4)--(z2)--(z5)--(z3);
\draw (z4)--(z3)--(z6);
\draw (y0)--(z0);
\draw (x0)--(q1);
\draw (y0)--(q2);
\draw (d1)--(d3);

\draw (d0) to[out=90,in=180, distance=0.5cm] (d1);
\draw (d0) to[out=270,in=180, distance=0.5cm] (x0);
\draw (d2) to[out=270,in=0, distance=0.5cm] (z0);
\draw (d2) to[out=90,in=0, distance=0.5cm] (d3);

\draw (x4) [fill=white] circle (\vr);
\draw (x1) [fill=black] circle (\vr);
\draw (x2) [fill=black] circle (\vr);
\draw (x3) [fill=black] circle (\vr);
\draw (x0) [fill=black] circle (\vr);
\draw (x5) [fill=black] circle (\vr);
\draw (x6) [fill=black] circle (\vr);
\draw (p1) [fill=black] circle (\vrs);
\draw (p2) [fill=black] circle (\vrs);
\draw (p3) [fill=black] circle (\vrs);
\draw (y4) [fill=white] circle (\vr);
\draw (y1) [fill=black] circle (\vr);
\draw (y2) [fill=black] circle (\vr);
\draw (y3) [fill=black] circle (\vr);
\draw (y0) [fill=black] circle (\vr);
\draw (y5) [fill=black] circle (\vr);
\draw (y6) [fill=black] circle (\vr);
\draw (z0) [fill=black] circle (\vr);
\draw (z4) [fill=white] circle (\vr);
\draw (z1) [fill=black] circle (\vr);
\draw (z2) [fill=black] circle (\vr);
\draw (z3) [fill=black] circle (\vr);
\draw (z6) [fill=black] circle (\vr);
\draw (z5) [fill=black] circle (\vr);

\end{tikzpicture}
\caption{A subcubic graph $G$ of order~$n$ satisfying $\ioc(G) = \frac{6}{7}n$. Dark vertices belong to the optimal code.}
\label{fig:cycle-halfgraphs}
\end{center}
\end{figure}

It is also interesting to notice that the above examples like half-graphs and the other graphs built using half-graphs, which serve as counterexamples to the bound in Theorem~\ref{thm:main1} for general graphs, contain induced $4$-cycles. Hence, even before proving a bound like the one in Theorem~\ref{thm:main1} in the most general case (where graphs contain $4$-cycles), one could also investigate if the same bound holds for graphs without \emph{induced} $4$-cycles. On the other hand, it would also be interesting to find graphs (if they exist) with non-induced $4$-cycles (that is, $4$-cycles with chords) which are counterexamples to the bound in Theorem~\ref{thm:main1}.

As a remark related to the application of the bound in Theorem~\ref{thm:main1} to other code numbers, it can be verified that any upper bound for IO-code numbers also holds for locating-total domination numbers of graphs. However, for the latter codes, there already exists a much stonger conjecture~\cite{FH-ElecJC-2016} in the case of twin-free graphs which states that the locating-total domination number of a twin-free graph on $n$ vertices is at most $\frac{2}{3}n$. Moreover, the best known approximation to the upper bound in the said conjecture is $\frac{3}{4}n$~\cite{FH-ElecJC-2016}. On the other hand, if one were to also allow twins in a graph $G$ of maximum degree at most $\Delta$, since any locating-total dominating set requires at least all-but-one twins in it, the $\frac{3}{4}n$-bound on the graph $G$ with its twins ``trimmed off'' and the inclusion of all-but-one twins in a locating-total dominating code roughly gives a $\left( \frac{\Delta-1}{\Delta}n \right)$-upper bound to a locating-total domination number of $G$. Notice that the latter is a better bound on a location-total domination number than a direct application of the bound in Theorem~\ref{thm:main1}.

We also recall here two questions by Henning and Yeo from~\cite{HeYe-14} about regular graphs. First, they conjectured that for every connected open twin-free cubic graph $G$, the upper bound $\ioc(G) \le \frac{3}{5}n$ holds, except for a finite number of graphs. Second, they asked whether $\ioc(G) \le \left(\frac{\Delta}{\Delta+1}\right)n$ holds for every connected open twin-free $\Delta$-regular graph $G$ of order $n$. 

\section{Acknowledgement}

The first two authors were supported by the French government IDEX-ISITE initiative CAP 20-25 (ANR-16-IDEX-0001), the International Research Center "Innovation Transportation and Production Systems" of the I-SITE CAP 20-25, and the ANR project GRALMECO (ANR-21-CE48-0004). The research of Dipayan Chakraborty was additionally supported in part by the Lebanese American University under the President's Intramural Research Fund PIRF0056. The research of Michael Henning was supported in part by the South African National Research Foundation, Grant Numbers 132588  and 129265, and the University of Johannesburg.

\end{document}